\documentclass[final,leqno,onefignum,onetabnum]{siamltex1213}

\usepackage{graphicx,epsfig}

\title{Bifurcation of limit cycles from a fold-fold singularity in planar switched systems\thanks{ This work was
supported by the Natural Science Foundation Grant
CMMI-1436856.}} 

\author{Oleg Makarenkov\thanks{Department of Mathematical Sciences, University of Texas at Dallas, 800 West Campbell Road 
Richardson, TX 75080 
(\email{makarenkov@utdallas.edu}). }}

\usepackage{hyperref}

\usepackage{ulem}

\usepackage{textcomp}

\usepackage[greek,english]{babel} 
\usepackage{amssymb,amsmath} 
\usepackage{amssymb,amsmath,color,latexsym} 
\usepackage{times}
\newcommand{\sign}{{\rm \hskip0.5pt sign \hskip1pt}}

\def\qed{\hfill $\square$}

\newcommand{\eps}{\varepsilon}

\newcommand{\R}{{\mathbb{R}}}

\usepackage{hyperref}

\newtheorem{remark}[theorem]{\rm \it Remark}

\begin{document}
\maketitle
\slugger{mms}{xxxx}{xx}{x}{x--x}

\begin{abstract}
We use a bifurcation approach to investigate the dynamics of a planar switched system that alternates between two smooth systems of ODEs denoted as (L) and (R) respectively. For an $x\in\R$ that plays the role of a bifurcation parameter, a switch to (R) occurs when the trajectory hits the switching line $\{x\}\times\R$ and a switch to (L) occurs when the trajectory hits the switching line $\{-x\}\times\R$. This type of switching is known as relay or hysteresis switching in control. The main result of the paper gives sufficient conditions for bifurcation of an attractive or repelling limit cycle from a point $O\in\{0\}\times\R$ when $x$ crosses 0. The result is achieved by spotting  a region where the dynamics of the system is described by the map $P(y)=y+\alpha y^3+\beta\dfrac{x}{y}+o(y).$ Motivated by applications to anti-lock braking systems, we focus on a particular class of switched systems where, for $x=0$, the point $O$ is a so-called fold-fold singularity, i.e. the vector fields of both systems (L) and (R) are parallel to $\{0\}\times \R$ at $O$. 
\end{abstract}

\begin{keywords} Switched system, relay system, fold-fold singularity, switched equilibrium, limit cycle, bifurcation, stability, normal form, anti-lock braking system\end{keywords}

\begin{AMS}37G15, 93C30\end{AMS}

\pagestyle{myheadings}
\thispagestyle{plain}
\markboth{Oleg Makarenkov}{\uppercase{Bifurcation of limit cycles from a fold-fold singularity}}

\section{Introduction}  
\begin{figure}[htbp]\label{f1}
\hskip0.15cm \includegraphics[scale=0.9]{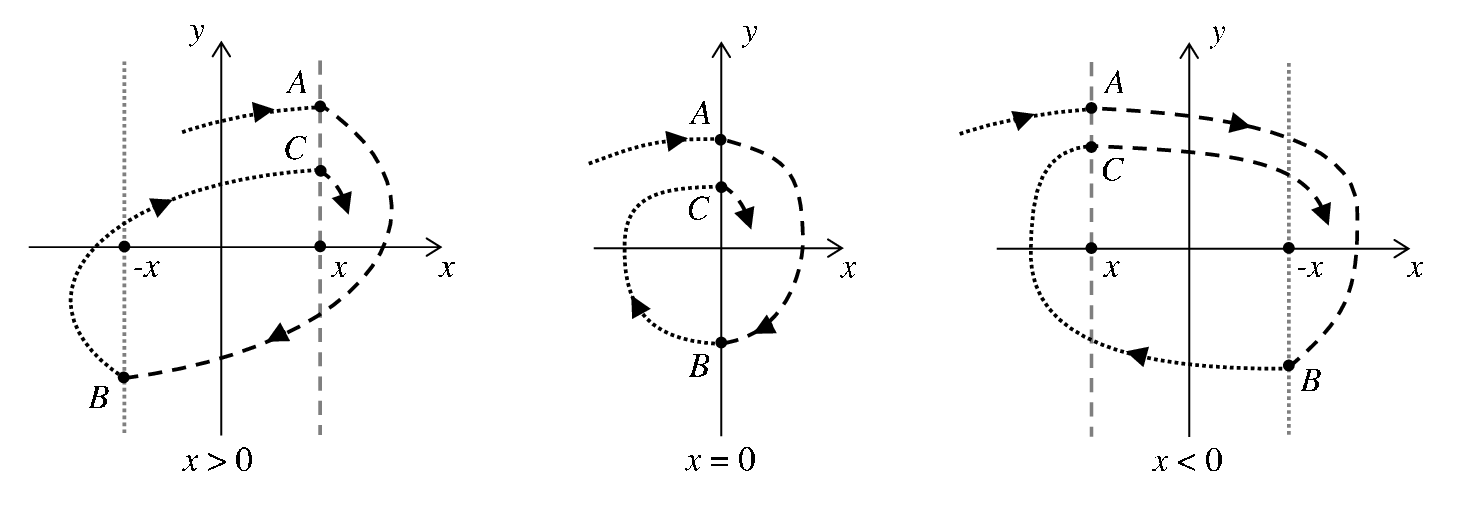}
\caption{Sample trajectory of system (\ref{fg})-(\ref{RL}) for different values of the parameter $x$. The dotted and dashed curves denote trajectories governed by system (\ref{fg}) with $i=L$ and $i=R$ respectively. The dotted line is where a switch to $i=L$ occurs. The dashed line is where a switch to $i=R$ occurs. The two lines coincide with $\{0\}\times\R$ when $x=0.$}
\end{figure}

This paper investigates the existence of attracting limit cycles in switched systems of the form
\begin{eqnarray}
&&  \begin{array}{l}
     \dot x(t)=f^i(x(t),y(t)),\\
    \dot y(t)=g^i(x(t),y(t)),
\end{array} \label{fg}\\
&&\begin{array}{l}
     i:=R,\quad{\rm if }\ \  x(t)=x,\\
     i:=L,\quad{\rm if }\ \ x(t)=-x,
  \end{array}\label{RL}
\end{eqnarray}
where $f^L,$ $f^R,$ $g^L,$ $g^R$ are smooth functions and $x\in\mathbb{R}$ is a parameter (we believe that using same letter $x$ to denote both the parameter and the solution doesn't cause a confusion and makes the notations less heavy). To draw a trajectory of system (\ref{fg})-(\ref{RL}) one  needs not just the initial point $(x(0),y(0))$, but also the letter of the system ($i=L$ or $i=R$) that governs the trajectory at $t=0.$ The trajectory is then governed by the system $i$ until it reaches one of the lines $\{-x\}\times\R$ or $\{x\}\times\R$, when $i$ switches to $i=L$ or $i=R$ according to whether $\{-x\}\times\R$ or $\{x\}\times\R$ is hit. The trajectory further travels along system $i$ until it reaches one of the switching lines again, when the same rule applies, see Figure~\ref{f1}.
In this paper we only deal with solutions that intersect the switching lines transversally.

\vskip0.2cm

\noindent The conditions for the existence of limit cycles in linear systems  (\ref{fg})-(\ref{RL}) are proposed in Astrom \cite{astrom} and Goncalves et al \cite{conc}. The case where (\ref{fg})-(\ref{RL}) admits a describing function is addressed in Tsypkin \cite{tsypkin}. These works found numerous applications in automatic control, where system (\ref{fg})-(\ref{RL}) is equivalently formulated as
$$
   \dot z=h(z,u)
$$
where $z$ is a vector and $u$ is a scalar control input that can only take a discrete set of values.
 The existence of a stable cycle for (\ref{fg})-(\ref{RL}) with piecewise linear $g^L$ and $g^R$ is established in Andronov et al \cite[Ch. III, \S5]{andr} in the context of clock modeling. Stable limit cycles in nonlinear systems of form (\ref{fg})-(\ref{RL}) are addressed in monographs by Barbashin \cite{barbashin} and Neimark \cite{neimark} along with applications in electromechanical engineering. Other applications where a  discrete-valued control is designed to produce limit cycles include switching power converters \cite{gupta}, intermittent therapy modeling in medicine \cite{tanaka}, grazing management in ecology \cite{meza}.  A particular model that motivated the current paper is an anti-lock braking system as introduced in \cite{tan09,pas06}. 
\vskip0.2cm

\noindent The paper deals with systems (\ref{fg})-(\ref{RL}) whose limit cycle shrinks to a point $(0,y_0)$ when $x\to 0.$ The latter can only happen when the vector fields $(f^L,g^L)$ and $(f^R,g^R)$ are oppositely directed at $(0,y_0)$, which is equivalent to saying that $(0,y_0)$ is an invisible equilibrium (or switched equilibrium) of the reduced system
\begin{equation}\label{np}
  \left(\begin{array}{c}
      \dot x(t)\\ \dot y(t)\end{array}\right)=\left\{\begin{array}{l}\left(\begin{array}{l}f^L(x(t),y(t))\\ g^L(x(t),y(t))\end{array}\right), \quad{\rm if}\ x(t)<0, \\
\left(\begin{array}{l}f^R(x(t),y(t))\\ g^R(x(t),y(t))\end{array}\right), \quad{\rm if}\ x(t)>0.\end{array}\right.
\end{equation}
We consider $y_0=0$ to shorten notations. When the vectors $(f^L(0),g^L(0))$ and $(f^R(0),g^R(0))$ are transversal to the switching manifold $\{0\}\times\R$, the bifurcation of a limit cycle of (\ref{fg})-(\ref{RL}) from an invisible equilibrium 0 of (\ref{np}) is of interest in power electronics. This type of bifurcations is a subject of a different paper. 
In this paper we address the situation where both $(f^L(0),g^L(0))$ and $(f^R(0),g^R(0))$ are parallel to the line $\{0\}\times\R$, which takes places in anti-lock braking systems and switched mechanical oscillators (see Sec.~\ref{appl}). By the other words, we consider the case where (\ref{np}) verifies
\begin{equation}\label{ff}
    f^L(0)=f^R(0)=0,
\end{equation}
i.e. where the origin is a fold-fold singularity of (\ref{np}). The main result of this paper is a sufficient condition that ensures bifurcation of a limit cycle of (\ref{fg})-(\ref{RL}) from a fold-fold singularity of (\ref{np}). This type of bifurcation is termed border-splitting bifurcation in \cite{lamb}.

\vskip0.2cm

\noindent System (\ref{fg})-(\ref{RL}) can be viewed as a result of a discontinuous perturbation of the switching manifold in system (\ref{np}). In this way, the current work complements the results by Guardia et al \cite{teixeira} and Kuznetsov et al \cite{kuznecov} on bifurcations of limit cycles from a fold-fold singularity of (\ref{np}) under smooth perturbations of the vector fields of (\ref{np}). Relevant results for continuous systems (\ref{np}) are obtained by Simpson-Meiss \cite{simpson} and  
Zou-Kuepper-Beyn \cite{zou}.

\vskip0.2cm

\noindent The paper is organized as follows. The next section is devoted to the proof of the main result, which is split in several steps. First, in section~\ref{ss1} we derive a normal form for the map that acts from the dashed line to the dotted line of Figure~\ref{f1} along the flow of an arbitrary smooth flow with a fold-singularity at the origin. This normal is obtained as a composition $\widetilde{\mathcal{P}}\circ\mathcal{P}$, where $\mathcal{P}$ is  the map from the dashed line to itself and $\widetilde{\mathcal{P}}$ is the rest of the trajectory $A\mapsto B$. The normal form for the map $\mathcal{P}$ turns out to be $\mathcal{P}(y)=-y+\alpha y^2$ in the neighborhood of the origin (lemma~\ref{mathcalP}) and the normal form for $\widetilde{\mathcal{P}}$ is found as $\widetilde{\mathcal{P}}(y)=y+\beta \dfrac{x}{y}$ in a special neighborhood of the origin (lemma~\ref{lem2}). This special neighborhood is colored in gray in Figure~\ref{fig2}. The results of section~\ref{ss1}  apply to the $R$-system of (\ref{fg}) in order to get the respective map $A\mapsto B$ and to the $L$-system of (\ref{fg}) in order to get the map $B\mapsto C.$ This allows to obtain (corollary~\ref{normalP}) the map $A\mapsto C$ as $P(y)=y+\alpha y^2+\beta\dfrac{x}{y}+o(y),$ whose dynamics is then investigated in section~\ref{dynamicsP}. Our main results theorem~\ref{mainthm} (section~\ref{dynamicsP}) then collects the findings of sections~\ref{ss1}-\ref{dynamicsP} in one statement about the limit cycles of the switched system (\ref{fg})-(\ref{RL}). Applications were a significant driving force towards this paper. A toy mechanical model that illustrates theorem~\ref{mainthm} is considered in section~\ref{appl1}, where following \cite{lamb} we just incorporated a relay (hysteresis) switching of the external force in the standard mass-spring oscillator, thus making the dynamics switching between two globally attracting equilibria. In section~\ref{appl2} we apply theorem~\ref{mainthm} to establish the existence of a limit cycle in an anti-lock braking system (ABS) with a two-phase controller of the pressure valves. Achieving such a limit cycle is a standard goal of an ABS controller as discussed in \cite{pas06}. The proofs of more technical statements (such as estimates for reminders) are sent to the appendix.

\vskip0.2cm

\noindent The proof of the main result (theorem~\ref{mainthm}) is split into a set of smaller statements (lemmas, corollaries and propositions). 
However, just formulations of those smaller statements will represent a proof of the main result for some readers. Indeed, spotting a suitable expansion and a suitable domain for the Poincar\'e map $P$  of the switched system (\ref{fg})-(\ref{RL}) (so that $P$ doesn't contain square-root terms) was the main difficulty in this paper. We stress that the approach of this paper doesn't permit us to analyze the dynamics outside the suitable domain mentioned above, where the effect of a square-root singularity may come into play (the interested reader can consult \cite{bristol} for the role of square-root maps and square-root singularities in the analysis of nonsmooth systems).

\section{The main result}\label{sec2}

The main result of this paper is achieved by analyzing the normal form of the Poincare map $y\mapsto P(y)$ of system  (\ref{fg})-(\ref{RL})  induced by the line $\{x\}\times\R$. To construct the map $P$, we consider the flow of system (\ref{fg}) with $i=R$ from the line $\{x\}\times\R$ to the line line $\{-x\}\times\R$. The respective map is denoted by $P^R_x$ and is called the point transformation. For example, $P^R_x(A)=B$ in the example of Figure~\ref{f1}. Furthermore, we denote by $P^L_{-x}$  the point transformation of (\ref{fg}) with $i=L$ from the line $\{-x\}\times\R$ to the line $\{x\}\times\R$. In particular, $P^L_{-x}(B)=C$ for the flow of Figure~\ref{f1}.
The Poincare map $P$ is obtained as a composition 
\begin{equation}\label{introduceP}
   P=P^L_{-x}\circ P^R_x,
\end{equation}
i.e. $P(A)=C$ for the points from Figure~\ref{f1}.
In section \ref{ss1} we derive normal forms for the point transformations $P^L_{-x}$ and $P^R_x$, that will allow us do draw the required conclusions about $P.$

\subsection{The normal form of  a point transformation in the neighborhood of a fold-fold singularity}\label{ss1}

\begin{figure}[h]
\begin{center}\includegraphics[scale=0.9]{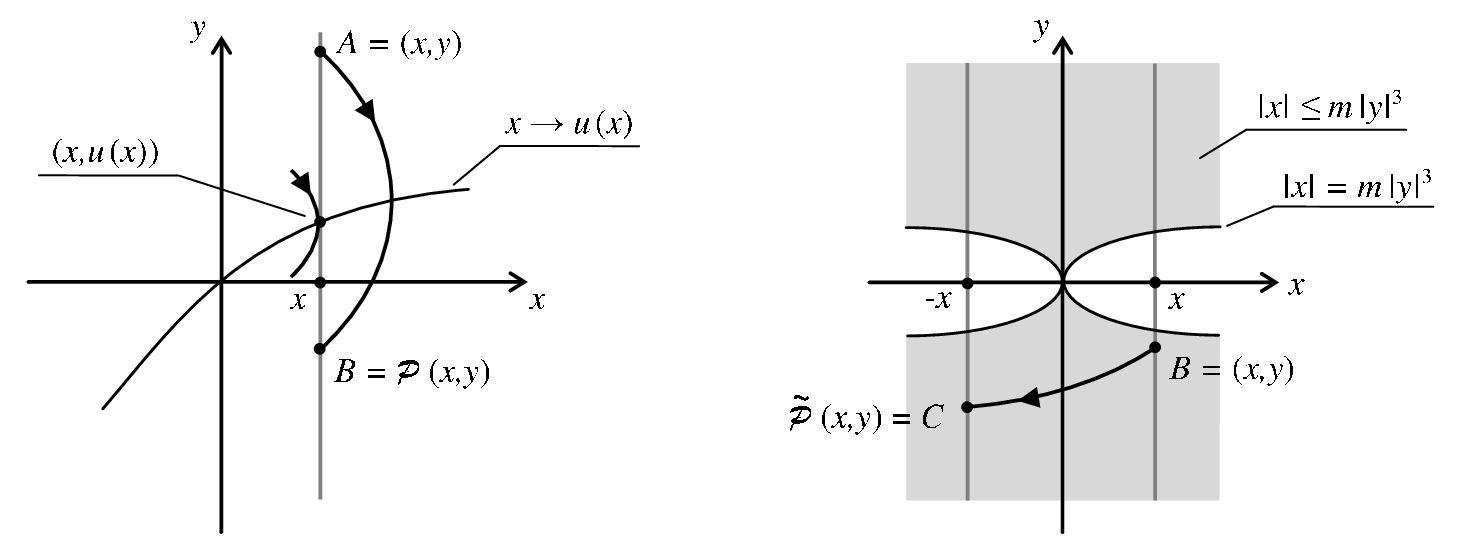}\end{center}
\caption{The two parts $A\mapsto B$ and $B\mapsto C$ of the point transformation $A\mapsto C$ of (\ref{neim}) induced by the lines $\{x\}\times\R$ and $\{-x\}\times\R$.}
  \label{fig2}
\end{figure}

Consider a planar system 
\begin{equation}\label{neim}
\begin{array}{l}
   \dot x=f(x,y),\\
   \dot y=g(x,y).
\end{array}
\end{equation}
Assuming that the origin is a fold-fold singularity of (\ref{neim}), we are going to spot the trajectories $t\mapsto x(t)$ of (\ref{neim}) that originate in $\{x\}\times\R$ at $t=0$ and reach $\{-x\}\times\R$ at $t=T$ while crossing $\{x\}\times\R$ at just one point $t_*\in (0,T).$  This will allow as to view the point transformation of (\ref{neim}) from a suitable subset of $\{x\}\times\R$ to $\{-x\}\times\R$ as a composition $\widetilde{\mathcal{P}}\circ\mathcal{P}$, where $\mathcal{P}$ and $\widetilde{\mathcal{P}}$ are defined as $\mathcal{P}(x(0))=x(t_*)$ and $\widetilde{\mathcal{P}}(x(t_*))=x(T)$ respectively, see Figure~\ref{fig2}.

\vskip0.2cm

\noindent The lemma~\ref{mathcalP}  below is saying that any solution of (\ref{neim}) which originates in $\{x\}\times\R$ reaches $\{x\}\times\R$ again in forward or backward time. This return is locally unique and the respective local map from $\{x\}\times\R$ to itself can be expanded as (\ref{expandable}).

\vskip0.2cm

\noindent In what follows, $t\mapsto (X(t,x,y),Y(t,x,y))$ denotes the general solution of (\ref{neim}) with the initial condition $(X(0,x,y),Y(0,x,y))=(x,y).$

\begin{lemma}\label{mathcalP} Assume that $C^4$ functions $f$ and $g$ satisfy
\begin{equation}\label{fdy}f(0)=0,\ \ f'_y(0)g(0)\not=0.\end{equation}
Let $u(x)\to 0$ as $x\to 0$ be the unique function such that $f(x,u(x))=0$ all $|x|\le\delta$, where $\delta>0$ is a suitable constant. The constant  $\delta>0$ can be diminished so that, for all $y\not=u(x)$, $|x|\le\delta$, and $|y|\le\delta$  the equation
\begin{equation}\label{solveX}
  X(t,x,y)=x
\end{equation}
admits a unique non-zero $C^2$
solution $T(x,y)\to 0$ as $y\to u(x).$ For this solution,
\begin{equation}\label{Tdif}
  T'_y(x,u(x))=-\dfrac{2}{g(x,u(x))}.
\end{equation}
The $C^2$ map 
$$
   {\mathcal P}(y)=Y(T(x,y),x,y)
$$ 
expands as
\begin{equation}\label{expandable}
  \mathcal{P}(y)= -y+\alpha y^2+\mathcal{R}(x,y),
\end{equation}
where 
\begin{equation}\label{R}
\alpha=2\dfrac{f'_x(0)+g'_y(0)}{g(0)}+\dfrac{f{}'_y{}'_y(0)}{f'_y(0)}, \ \ \mathcal{R}(0,0)=\mathcal{R}'_x(0,0)=\mathcal{R}'_y(0,0)=\mathcal{R}'_y{}'_y(0,0)=0.
\end{equation}
\end{lemma}

{\bf \it Proof.} {\it Step 1:} {\it The existence of $T(x,y)$.} 
 The existence of $u(t)$ under condition (\ref{fdy}) follows from the Implicit Function Theorem (see \cite[\S~8.5.4, Theorem~1]{zorich}). To solve (\ref{solveX}) we expand $X$ in 
Taylor series as
$$
  X(t,x,y)=x+X'_t(0,x,y)t+X'_t{}'_t(0,x,y)\dfrac{t^2}{2}+\rho(t,x,y)t.
$$
Equation (\ref{solveX}) takes the form
\begin{equation}\label{Req}
   X'_t(0,x,y)+X'_t{}'_t(0,x,y)\dfrac{t}{2}+\rho(t,x,y)=0,
\end{equation}
where (see Appendix~\ref{rhoprop}) uniformly in $(x,y)\in[-\delta,\delta]^2$ 
\begin{eqnarray*}
&&
  \lim_{t\to 0} \rho(t,x,y)=\lim_{t\to 0} \rho'_t(t,x,y)=\lim_{t\to 0} \rho'_y(t,x,y)= \lim_{t\to 0} \rho'_t{}'_y(t,x,y)=\lim_{t\to 0} \rho'_y{}'_y(t,x,y)=0,\\  
&&\lim_{t\to 0} \rho'_t{}'_t(t,x,y)=\dfrac{x'_t{}'_t{}'_t(0,x,y)}{3},\quad \lim\limits_{t\to 0}
\rho'_x(t,x,y)=\lim\limits_{t\to 0}\rho{}'_t{}'_x(t,x,y)=\lim\limits_{t\to 0}\rho{}'_x{}'_x(t,x,y)=0.
\end{eqnarray*}
Since $X{}'_t{}'_t(0,0,0)=f{}'_y(0)g(0),$ see Appendix~\ref{partial}, 
by the Implicit Function Theorem there exists a unique $T(x,y)\to 0$ as $(x,y)\to 0$ that solves (\ref{solveX}) in $[-\delta,\delta]^2.$ 
Since $T(x,u(x))=0$, then  $T(x,y)\not=0$ for all $y\in[-\delta,u(x))\cup(u(x),\delta]$.
Replacing $t$ by $T(x,y)$ in (\ref{Req}) and by taking the derivative with respect to $x$ and $y$ one can express $T'_x(x,y)$ and $T'_y(x,y)$ as
\begin{eqnarray*}
  T'_x(x,y)&=&\dfrac{X'_t{}'_x(0,x,y)+(1/2)X'_t{}'_t{}'_x(0,x,y)T(x,y)+\rho'_x(T(x,y),x,y)}{-(1/2)X'_t{}'_t(0,x,y)-\rho'_t(T(x,y),x,y)},\\
  T'_y(x,y)&=&\dfrac{X'_t{}'_y(0,x,y)+(1/2)X'_t{}'_t{}'_y(0,x,y)T(x,y)+\rho'_y(T(x,y),x,y)}{-(1/2)X'_t{}'_t(0,x,y)-\rho'_t(T(x,y),x,y)},
\end{eqnarray*}
from where  $T$ is $C^2.$ Moreover, using the formulas from Appendix~\ref{partial}, we conclude that 
 $T'_y(x,u(x))$ is given by (\ref{Tdif}) and that  
 $T{}'_y{}'_y(0)$ is given by
$$   T{}'_y{}'_y(0)=4\dfrac{f'_x(0)+g'_y(0)}{g(0)^2}+2\dfrac{f{}'_y{}'_y(0)}{f'_y(0)g(0)}.
$$

\noindent {\it Step 2:} {\it The expansion of $\mathcal{P}(y)$.} From $T(x,u(x))=0$ and $Y'_t(0,x,u(x))T'_y(x,u(x))=-2$ one has 
$$\mathcal{P}(u(x))=u(x)\quad {\rm and}\quad  
 \mathcal{P}'(u(x))=-1.$$ Computing $\mathcal{P}''$ yields
\begin{equation}\label{0fo}
\begin{array}{lll}
\dfrac{1}{2}\mathcal{P}''(0)&=&\dfrac{1}{2}\left(Y{}'_t{}'_t(0)T'_y(0)^2+Y'_t(0)T{}'_y{}'_y(0)+2Y{}'_t{}'_y(0)T'_y(0)+Y{}'_y{}'_y(0)\right)=\\
&&=\dfrac{1}{2}g(0)T{}'_y{}'_y(0)=\alpha.
\end{array}
\end{equation}
Therefore, the Taylor series for $\mathcal{P}(y)$ around $y=u(x)$ is
\begin{equation}\label{1fo}
\begin{array}{rcl}
  \mathcal{P}(y)&=&\mathcal{P}(u(x))+\mathcal{P}'(u(x))(y-u(x))+\dfrac{1}{2}\mathcal{P}''(y_*)(y-u(x))^2=\\
& =&2u(x)-y+\dfrac{1}{2}\mathcal{P}''(y_*)(y-u(x))^2\ =:\ -y+\alpha y^2+\mathcal{R}(x,y).
\end{array}
\end{equation}
The $C^2$ smoothness of $\mathcal{R}$ follows from the $C^2$ smoothness of $\mathcal{P}$, which is a consequence of $C^2$ smoothness of $T$. Properties (\ref{R}) follow by direct computation executed for $\mathcal{R}(x,y)=\mathcal{P}(y)+y-\dfrac{1}{2}\mathcal{P}''(0)y^2.$\qed

\vskip0.2cm

\noindent The next lemma is saying that any solution of (\ref{neim}) that originates at $(x,y)$ with $|x|\le m|y|^3$ reaches $\{-x\}\times\R$ in forward or backward time. This return is locally unique and the respective local map from $\{x\}\times\R$ to $\{-x\}\times\R$ can be expanded on $|x|\le m|y|^3$ as (\ref{widetildeP}).

\vskip0.2cm

\begin{lemma}\label{lem2} Assume that both $f$ and $g$ are $C^4$ and $$f(0)=0,\ \ f{}'_y(0)\not=0.$$ Then, for any $m>0$ there exists $\delta>0$ such that 
for all $|x|\le m|y|^3$ and $|y|\le \delta$ 
the equation
\begin{equation}\label{eqX}
   X(\widetilde T,x,y)=-x
\end{equation}
admits a unique solution 
\begin{equation}\label{uniquesolution}|\widetilde T(x,y)|\le \dfrac{4m}{|f'_y(0)|}y^2.\end{equation}
 This solution is $C^2$ in $0<\sqrt[3]{|x|/m}<|y|<\delta$ and   
\begin{equation}\label{Tas}
   \widetilde T(x,y)\cdot\dfrac{y}{x}\to -\dfrac{2}{f'_y(0)}\quad {\rm as}\quad |x|\le m|y|^3,\ y\to 0.
\end{equation}
The respective $C^2$ map
$$
  \widetilde {\mathcal{P}}(y)=Y(\widetilde{T}(x,y),x,y)
$$
expands as
\begin{equation}\label{widetildeP}
   \widetilde {\mathcal{P}}(y)=y+\beta \dfrac{x}{y}+\widetilde{\mathcal R}(x,y),
\end{equation}
where
\begin{equation}\label{tildeR}
\begin{array}{l}
    \lim\limits_{|x|\le m |y|^3,\ y\to 0}\dfrac{\widetilde {\mathcal R}(x,y)}{y^2}=0,\ \beta=-\dfrac{2g(0)}{f'_y(0)}\\ \lim\limits_{|x|\le m |y|^3,\ y\to 0}\dfrac{\widetilde {\mathcal{R}}'_y(x,y)}{y}=0,\ \lim\limits_{|x|\le m |y|^3,\ y\to 0}\widetilde {\mathcal{R}}'_x(x,y)y=0.
\end{array}
\end{equation}
Furthermore,
\begin{equation}\label{A2}
\begin{array}{l}
   \widetilde{\mathcal{P}}\left(\left[\sqrt[3]{|x/m|},\delta\right]\right)\subset \left[(1-\eps)\sqrt[3]{|x/m|},(1+\eps)\delta\right],\\ \widetilde{\mathcal{P}}\left(\left[-\delta,-\sqrt[3]{|x/m|}\right]\right)\subset \left[-(1+\eps)\delta,-(1-\eps)\sqrt[3]{|x/m|}\right],
\end{array}\quad{for\ all\ }|x|\le m\delta^3. 
\end{equation}
\end{lemma}

\vskip0.2cm

{\it Proof.} {\it Step 1: The existence and uniqueness of $\widetilde T$}. Introduce
$$
  F(t,x,y)=\dfrac{X(t,x,y)+x}{y}.
$$
Then $|F(0,x,y)|=|2x/y|\le 2my^2$ when $|x|\le m |y|^3$. Furthermore, since
$$
   \dfrac{X{}'_t(0,x,y)}{y}\to X{}'_t{}'_y(0,x,y)\quad{\rm as}\quad |x|\le m|y|^3,\ \ y\to 0
$$ 
then, see also Appendix~\ref{partial},
 $F'_t(t,x,y)\to X'_t{}'_y(0)=f'_y(0)$ as $|t|\le k |y|^2,$ $|x|\le m |y|^3,$ $y\to 0.$ The latter holds for any fixed $k>0$ and, in particular, for $k=4m/f'_y(0)$. Thus, the existence, uniqueness, differentiability and the estimate for $\widetilde T$ follow from Theorem~\ref{implicit} applied with $n=2,$ $\xi=t$, $\zeta=x,$ $M=2m$, and $q=\dfrac{1}{2}.$

\vskip0.2cm

\noindent {\it Step 2: The asymptotic of $\widetilde T(x,y)x/y.$}   Expanding $X(t,x,y)=x+X'_t(t_*,x,y)t$, we can rewrite (\ref{eqX}) as
$$
   2x+X'_t(t_*,x,y)\widetilde{T}(x,y)=0,
$$
from where 
$$
  \widetilde T(x,y)\cdot\dfrac{y}{x}=-\dfrac{2}{\dfrac{X'_t(t_*,x,y)}{y}}=-\dfrac{2}{X{}'_t{}'_t(t_{**},x,y)\dfrac{t_*}{y}+\dfrac{X'_t(0,x,y)}{y}}.
$$

\vskip0.2cm

\noindent {\it Step 3: The asymptotic of $\widetilde {\mathcal {R}}(x,y)y$.} 
Expanding $X(t,x,y)=x+\bar\Delta(t,x,y)$, we can rewrite (\ref{eqX}) as
$   2x+\bar \Delta(\widetilde{T}(x,y),x,y)=0,
$ from where
$$
   \widetilde T'_x(x,y)y=\dfrac{-2-\bar\Delta'_x(\widetilde T(x,y),x,y)}{\dfrac{\bar\Delta'_t(\widetilde T(x,y),x,y)}{y}}\to-\dfrac{2}{X{}'_t{}'_y(0)}\quad{\rm as}\quad |x|\le m|y|^3,\ \ y\to 0.
$$
Therefore,
$$
\begin{array}{l}
   \widetilde {\mathcal{R}}'_x(x,y)y=\left[Y'_t(\widetilde T(x,y),x,y)\widetilde T'_x(x,y)+Y'_x(\widetilde T(x,y),x,y)-\beta\dfrac{1}{y}\right]y\to\\
\to Y'_t(0)\left(-\dfrac{2}{X{}'_t{}'_y(0)}\right)-\beta=0\quad{\rm as}\quad |x|\le m|y|^3,\ \ y\to 0.
\end{array}
$$

\vskip0.2cm

\noindent {\it Step 4: The asymptotics of $\dfrac{\mathcal{R}(x,y)}{y^2}$ and $\dfrac{\mathcal{R}'_y(x,y)}{y}$.} Expanding $X(t,x,y)=x+X'_t(0,x,y)t+\widetilde\Delta(t,x,y)$, we can rewrite (\ref{eqX}) as
\begin{equation}\label{2x}
   2x+X'_t(0,x,y)\widetilde{T}(x,y)+\widetilde \Delta(\widetilde{T}(x,y),x,y)=0,
\end{equation}
and so
$$
   \widetilde{T}(x,y)=-\dfrac{2x}{X'_t(0,x,y)}-\dfrac{\widetilde \Delta(\widetilde{T}(x,y),x,y)}{X'_t(0,x,y)}.
$$
By using
\begin{equation}\label{thfor}
  \widetilde{\mathcal{P}}(x,y)=y+Y'_t(t_*,x,y)\widetilde T(x,y)=y-\dfrac{2Y'_t(0)}{X{}'_t{}'_y(0)}\cdot\dfrac{x}{y}+\widetilde {\mathcal{R}}(x,y)
\end{equation}
one now obtains the following formula for $\widetilde {\mathcal R}(x,y)$
$$
  \widetilde {\mathcal{R}}(x,y)=-\dfrac{2Y'_t(t_*,x,y)x}{X'_t(0,x,y)}-\dfrac{Y'_t(t_*,x,y)\widetilde \Delta(\widetilde T(x,y),x,y)}{X'_t(0,x,y)}+\dfrac{2Y'_t(0)}{X{}'_t{}'_y(0)}\cdot\dfrac{x}{y},
$$
from where the required limiting relations (\ref{tildeR}) follow, see Appendix~\ref{appendixC}.

\noindent {\it Step 5: The relations (\ref{A2}).}
Let us $\delta>0$ be so small that $\dfrac{|\beta x/y+\widetilde{\mathcal{R}}(x,y)|}{|y|}\le\eps$ for all $x\le m|y|^3$ and $|y|\le\delta.$ For these values of $x$ and $y$, the expansion (\ref{widetildeP}) yields
$$
y(1-\eps\cdot{\rm sign}|y|)\le\widetilde{\mathcal{P}}(y)\le y(1+\eps\cdot{\rm sign}|y|),
$$
which implies the required inclusions for the values of $\widetilde{\mathcal{P}}.$

 \qed

\vskip0.2cm

\noindent Lemma~\ref{lem2} provides information about those trajectories of (\ref{neim}) only whose initial conditions $(x,y)$ satisfy $|x|\le m |y|^3$. We will need similar properties of the reminder in lemma~\ref{mathcalP} when studying the composition $\widetilde{\mathcal{P}}\circ\mathcal{P}$. These properties are given by the following corollary.

\vskip0.2cm

\begin{corollary}\label{cor2_3} Under the conditions of lemma~\ref{mathcalP}, for any $m>0$,
\begin{equation}\label{SS}
\begin{array}{l}
\lim\limits_{{|x|}\le m|y|^3,\ y\to 0}\dfrac{T(x,y)}{y}=-\dfrac{2}{g(0)},\\
\lim\limits_{{|x|}\le m|y|^3,\ y\to 0}\dfrac{\mathcal{R}(x,y)}{y^2}=0,\quad \lim\limits_{{|x|}\le m|y|^3,\ y\to 0}\dfrac{{\mathcal{R}}'_y(x,y)}{y}=0.
\end{array}
\end{equation}
In particular,  for any  $\eps>0$ there exists $\delta>0$ such that  
\begin{equation}\label{A1}
\begin{array}{l}
{\mathcal P}\left(\left[\sqrt[3]{|x/m|},\delta\right]\right)\subset \left[-(1+\eps)\delta,-(1-\eps)\sqrt[3]{|x/m|}\right], \\
{\mathcal P}\left(\left[-\delta,-\sqrt[3]{|x/m|}\right]\right)\subset \left[(1-\eps)\sqrt[3]{|x/m|},(1+\eps) \delta\right],
\end{array}\quad for\  all \quad  |x|\le m\delta^3. 
\end{equation}
\end{corollary}

\vskip0.2cm

{\it Proof.} The proof of (\ref{SS}) is sent to Appendix~\ref{reminderR}. Here we just focus on the relations (\ref{A1}).
Let us $\delta>0$ be so small that $\dfrac{|\alpha y^2+\mathcal{R}(x,y)|}{|y|}\le\eps$ for all $x\le m|y|^3$ and $|y|\le\delta.$ For these values of $x$ and $y$, the expansion (\ref{expandable}) yields
$$
-y(1+\eps\cdot{\rm sign}|y|)\le\mathcal{P}(y)\le -y(1-\eps\cdot{\rm sign}|y|),
$$
which implies the required inclusions for the values of $\mathcal{P}.$ \qed

\vskip0.2cm

\noindent We are now in the position to combine the maps $\mathcal{P}$ and $\widetilde{\mathcal{P}}$ and to relate the composition $\widetilde{\mathcal{P}}\circ\mathcal{P}$ to the flow of system (\ref{neim}).

\vskip0.2cm

\begin{corollary}\label{compo} Assume that the conditions of lemmas \ref{mathcalP} and \ref{lem2} hold and let $\mathcal{P}$ and $\widetilde{\mathcal{P}}$ be the maps provided by these lemmas. Then 
$$
    \widetilde{\mathcal{P}}(\mathcal{P}(y))=-y+\alpha y^2-\beta\dfrac{x}{y}+r(x,y),
$$
where, for any $m>0$, the map $r$ is $C^2$ in $0<\sqrt[3]{|x/m|}<|y|$ and 
\begin{equation}\label{limrel}
\lim\limits_{|x|\le m |y|^3,\ y\to 0} r(x,y)= \lim\limits_{|x|\le m |y|^3,\ y\to 0}\dfrac{r'_y(x,y)}{y}=\lim\limits_{|x|\le m |y|^3,\ y\to 0}r'_x(x,y)y=0.
\end{equation}
For any $m>0$ there exists $\delta>0$ such that, for 
${|x|}<m\delta^3$ and 
\begin{equation}\label{signx}
{\rm sign}(x)=-{\rm sign}(f'_y(0)g(0)),
\end{equation}
the map 
$ \widetilde{\mathcal{P}}\circ\mathcal{P}$ describes the transformation of the interval \begin{equation}\label{signy}
\{x\}\times[-\sign(g(0))\sqrt[3]{|x/m|},-\sign(g(0))\delta]
\end{equation} 
of the line $\{x\}\times\R$ to the line $\{-x\}\times\R$ under the action of the flow of (\ref{neim}). In particular, $T(x,y)>0$ and $\widetilde T(-x,\mathcal{P}(x,y))>0$ for $(x,y)$ from (\ref{signx})-(\ref{signy}).
\end{corollary}

\vskip0.2cm

{\it Proof.} {\it Step 1: Properties of $r$.} Direct computation leads to
\begin{eqnarray*}
&& r(x,y)=\widetilde{\mathcal{R}}\left(x,-y+\alpha y^2+\mathcal{R}(x,y)\right)+\mathcal{R}(x,y)+\widetilde{\mathcal{R}}(x,y)+\beta\dfrac{x}{y}+\beta\dfrac{x}{-y+\alpha y^2+\mathcal{R}(x,y)}
\end{eqnarray*}
and the required properties of $r$ follow  from the respective properties (\ref{R}), (\ref{SS}) and (\ref{tildeR}) of $\mathcal{R}$ and $\widetilde{\mathcal{R}},$ see Appendix~\ref{reminderr}. 

\vskip0.2cm

\noindent {\it Step 2: The relation between the map $\widetilde{\mathcal{P}}\circ\mathcal{P}$ and the flow of (\ref{neim}).} The map $\widetilde{\mathcal{P}}\circ\mathcal{P}$ is a point transformation for the flow of (\ref{neim}), if the following two properties hold.
\begin{itemize}
\item[1)] Both $T(x,y)$ and $\widetilde T(x,y)$ that appear in the definitions of 
$\widetilde{\mathcal{P}}$ and $\mathcal{P}$ in lemmas 
\ref{mathcalP} and \ref{lem2} are positive. 

{\it Proof.} From (\ref{Tdif}) we have $\sign(T(x,y))=-\sign(g(0)y)$, i.e. $\widetilde T(x,y)>0$, iff $\sign(y)=-\sign(g(0))$ regardless of the value of $x.$ The latter leads to (\ref{signy}). From  
(\ref{Tas}) we conclude that $\sign(\widetilde T(x,y))=-\sign(f'_y(0)xy)$, i.e. $\widetilde T(x,P(y))>0$, iff $\sign(x)=-\sign(f'_y(0)P(y))$, which gives (\ref{signx}).
\item[2)] $X(t,x,y)\not=-x$ when $t\in(0,T(x,y)+\widetilde T(x,y)).$

{\it Proof.} $X(t,x,y)\not=-x$ for all $t\in[T(x,y),\widetilde T(x,y))$ by the uniqueness of $\widetilde T(x,y)$ ensured by lemma~\ref{lem2}.
It is therefore sufficient to prove that $X(t,x,y)\not=-x$ when $t\in(0,T(x,y))$. To have the latter it is sufficient to check that the curve $\cup_{t\in(0,T(x,y))}\{(X(t,x,y),Y(t,x,y))\}$ and the line ${-x}\times\R$ are located on different sides of the line ${x}\times\R$. This property holds, if 
\begin{equation}\label{SSSS}
\sign(x)=-\sign(f(x,y)), \quad{\rm for}\ (x,y)\ {\rm from}\ (\ref{signy}).
\end{equation} For small $|y|>0$, property (\ref{SSSS}) takes the form $\sign x=\sign (f'(0)y)$, which does indeed hold when $(x,y)$ satisfies
(\ref{signx})-(\ref{signy}).
\end{itemize}

\qed

\subsection{The normal form of the Poincare map in the neighborhood of a fold-fold singularity}\label{sec22}

\noindent Let $P^R_x$ be the composition $\widetilde{\mathcal{P}}\circ\mathcal{P}$ obtained by applying lemmas \ref{mathcalP} and \ref{lem2} to (\ref{fg}) with $i=R.$ Let $P^L_{-x}$ be  the composition $\widetilde{\mathcal{P}}\circ\mathcal{P}$ obtained by applying lemmas \ref{mathcalP} and \ref{lem2} to (\ref{fg}) with $i=L$ and by further replacing $x$ by $-x$. Using so defined $P_{-x}^L$ and $P_x^R$, introduce $P$ according to (\ref{introduceP}).

\vskip0.2cm

\begin{corollary}\label{normalP} Assume that both 
(\ref{fg}) with $i=L$ and (\ref{fg}) with $i=R$ satisfy the assumptions of lemmas \ref{mathcalP} and \ref{lem2}. Then,
for any $m>0$ there exists $\delta>0$ such that 
for all $|x|\le m|y|^3$ and $|y|\le \delta$ 
the map (\ref{introduceP}) admits a representation $$P(y)=y+(\alpha^L-\alpha^R)y^2+(\beta^R-\beta^L)\dfrac{x}{y}+\Delta(x,y),$$
where
$$
\lim\limits_{|x|\le m |y|^3,\ y\to 0}\Delta(x,y)= \lim\limits_{|x|\le m |y|^3,\ y\to 0}\dfrac{\Delta'_y(x,y)}{y}=\lim\limits_{|x|\le m |y|^3,\ y\to 0}\Delta'_x(x,y)y=0.
$$
If, in addition, 
\begin{equation}\label{inadd}
   \sign(x)=-\sign(f^R{}'_y(0)g^R(0)),\quad f^R{}'_y(0)f^L{}'_y(0)>0,\quad g^R(0)g^L(0)<0, 
\end{equation}
then for 
\begin{equation}\label{thenfor}
y\in \left\{y:-\sign(g^R(0))y\in\left[\sqrt[3]{|x/m|},\delta\right]\right\}\end{equation} 
the
map $y\mapsto P(y)$ 
is the Poincar\'e map for  switched system
(\ref{fg})-(\ref{RL}) induced by the cross-section $\{x\}\times\R.$
\end{corollary}

\vskip0.2cm

{\it Proof.} {\it Step 1: Properties of $\Delta(x,y)$.} Corollary~\ref{compo} implies that
\begin{eqnarray*}
P^R_x(y)&=&-y+\alpha^R y^2-\beta^R \dfrac{x}{y}+\Delta^R(x,y),\\
P^L_{-x}(y)&=&-y+\alpha^L y^2+\beta^L \dfrac{x}{y}+\Delta^L(x,y),
\end{eqnarray*}
where both $r(x,y)=\Delta^R(x,y)$ and $r(x,y)=\Delta^L(x,y)$ satisfy the limiting relations (\ref{limrel}).
Direct computation yields
\begin{eqnarray*}
   \Delta(x,y)&=&\Delta^L(x,y)-\Delta^R(x,y)+\alpha^L \Delta^R(x,y)^2+2\alpha^L\beta^R x+\alpha^L(\beta^R)^2\dfrac{x^2}{y^2}-\\
& &-2\alpha^L\beta^R \Delta^R(x,y)\dfrac{x}{y} 
-2\alpha^L\Delta^R(x,y)y-2\alpha^L\alpha^R\beta^Rxy+2\alpha^L\alpha^R\Delta^R(x,y)y^2-\\
&&- 2\alpha^L\alpha^R y^3+\alpha^L(\alpha^R)^2y^4+\dfrac{\beta^L x}{\Delta^R(x,y)-\beta^R(x/y)-y+\alpha^R y^2}+\beta^L\dfrac{x}{y}+\\
&&+\Delta^L\left(x,-y+\alpha^R y^2-\beta^R\dfrac{x}{y}+\Delta^R(x,y)\right)
\end{eqnarray*}
and the required properties of $\Delta$ follow from (\ref{limrel}) along the lines of the proofs carried out in Appendix~\ref{reminderr}.

\vskip0.2cm

\noindent {\it Step 2: The relation between the map $P$ and the flow of switched system
(\ref{fg})-(\ref{RL}).} According to corollary~\ref{compo} the map $P_x^R$ is a point transformation for (\ref{fg}) with $i=R$ on (\ref{thenfor}), where $\sign(x)=-\sign\left((f^R )'_y(0)g^R(0)\right).$

\vskip0.2cm

\noindent The map $P$ is a Poincar\'e map for (\ref{fg})-(\ref{RL}), if the following two properties hold.
\begin{itemize}
\item[1)]  $P^L_{-x}$ is a point transformation for (\ref{fg}) with $i=L$ from a subset $W$ of $\{-x\}\times\R$ to $\{x\}\times\R$. 

\vskip0.2cm

{\it Proof.} Based on corollary~\ref{compo}, it suffices to have $\sign(-x)=-\sign\left((f^L)'_y(0)g^L(0)\right)$ for 1) to hold. This sufficient property follows from (\ref{inadd}). Moreover, corollary~\ref{compo} implies that we can take 
$W=\left\{y:-\sign(g^L(0))y\in[\sqrt[3]{|x/M|},\delta L]\right\}$ where $M>0$ and $L>0$ can be chosen arbitrary.

\item[2)] $W\supset  W_0=P_x^R\left(\left\{y:-\sign(g^R(0))y\in\left[\sqrt[3]{|x/m|},\delta\right]\right\}\right).$

{\it Proof.} By (\ref{A1}) and (\ref{A2}) for any $\eps>0$ there exists $\delta>0$ such that
\begin{equation}\label{SQ1}
   W_0\subset \left\{y:\sign(g^R(0))y\in\left[(1-\eps)\sqrt[3]{|x/m|},(1+\eps)\delta\right]\right\},
\end{equation}
for all $|x|\le m|y|^3,$ $|y|\le\delta.$ Select $M\in(0,m).$ Then
\begin{equation}\label{SQ2}
   \left[(1-\eps)\sqrt[3]{|x/m|},(1+\eps)\delta\right]\subset\left[\sqrt[3]{|x/M|},(1+\eps)\delta\right]
\end{equation}
for all $|x|\le m\delta^3$ provided that $\delta>0$ is sufficiently small. Since, by (\ref{inadd}),  $g^R(0)g^L(0)<0,$ the property (\ref{SQ1})-(\ref{SQ2}) imply $W\supset W_0.$
\end{itemize}

\qed

\begin{remark}\label{note}{\rm Note that since the cross-section is selected as $\{x\}\times\mathbb{R}$, then the solutions of (\ref{fg})-(\ref{RL}) under consideration are governed by (\ref{fg}) with $i=R$ at the initial time. In particular, there is no uncertainty with the choice of the initial letter.}
\end{remark}

\subsection{The dynamics of the map $\mathbf{P(y)=y+\alpha y^2+\beta \dfrac{x}{y}+o(y)}$}\label{dynamicsP}

$\ $

When the term $x/y$ is replaced by $x$, the map $P$ turns into the normal form of a so-called fold bifurcation, see Kuznetsov \cite[\S~3.2]{kuznecov1}, whose dynamics is addressed in the literature starting from Guckenheimer \cite{guck}.

\noindent \begin{proposition}\label{dynamicsP} Consider a map
$$
  P(y)=y+\alpha y^2+\beta \dfrac{x}{y}+R(x,y),
$$
such that, for any $m>0$, the $C^2$ in $0<\sqrt[3]{|x/m|}<|y|$ reminder $R$ verifies
$$
\lim\limits_{|x|\le m |y|^3,\ y\to 0}R(x,y)= \lim\limits_{|x|\le m |y|^3,\ y\to 0}\dfrac{R'_y(x,y)}{y}=\lim\limits_{|x|\le m |y|^3,\ y\to 0}R'_x(x,y)y=0.
$$
If $\alpha\beta\not=0$, then for any $m>\dfrac{|\alpha|}{|\beta|}$ there exist $\delta>0$ and $\gamma>0$ such that for all $|x|<\gamma$ the map $P$ admits a unique fixed point $y(x)$ in the set
$$ \left[-\delta,-\sqrt[3]{|x|/m}\right]\cup\left[\sqrt[3]{|x|/m},\delta\right]=:I_{-1}\cup I_{1}.$$
Moreover,
\begin{equation}\label{moreover}
   \dfrac{x}{y(x)^3}\to -\dfrac{\alpha}{\beta}\quad{\rm as}\quad x\to 0,
\end{equation}	
in particular,
\begin{equation}\label{ie}
y(x)\in I_{{\rm sign}(-\alpha\beta x)}.
\end{equation}
If $x\beta>0$, then $y(x)$ is an attractor whose domain of attraction is at least $I_{{\rm sign}(-\alpha\beta x)}$. If $x\beta<0$, then $y(x)$ is a repeller and each trajectory that originates in $I_{{\rm sign}(-\alpha\beta x)}\backslash\{y(x)\}$ leaves $I_{{\rm sign}(-\alpha\beta x)}$ in finite time. 
\end{proposition}

\vskip0.2cm

{\it Proof.} {\it Step 1. The existence and uniqueness of $x(y)$.} The equation $P(y)=y$ is equivalent to $F(x,y)=0,$ where
$$
  F(x,y)=\alpha y^3+\beta x+R(x,y)y,
$$
so that $F(0,0)=0.$ We will use theorem~\ref{implicit} to prove the existence and uniqueness of zeros of $F(x,y).$ We have 
$$
    F'_x(x,y)\to\beta\quad{\rm as}\quad |x|\le m|y|^3,\ y\to 0.
$$
Furthermore, for any $\sigma>0$ there exists $\delta>0$ such that
\begin{equation}\label{the0}
   |F(0,y)|=\left(|\alpha|+\left|\dfrac{R(x,y)}{y^2}\right|\right)|y|^3\le\left(|\alpha|+\sigma\right)|y|^3,\quad|y|\le\delta.
\end{equation}
To fulfill (\ref{M})-(\ref{Lq}) it is, therefore, sufficient to find $\sigma>0$ and $q\in(0,1)$ such that 
$$
  \left(|\alpha|+\sigma\right)\le |\beta|mq.
$$
By the conditions of the proposition, $m=\dfrac{|\alpha|}{|\beta|}+r,$ where $r>0.$
Therefore,
\begin{equation}\label{the}
  \sigma\le|\beta|mq-|\alpha|=|\beta|\left(\dfrac{|\alpha|}{|\beta|}+r\right)q-|\alpha|=|\alpha|q+|\beta|rq-|\alpha|=-|\alpha|(q-1)+|\beta|rq,
\end{equation}
whose right-hand side can be made strictly positive by choosing $q\in(0,1)$ sufficiently close to 1. Let $\sigma>0$ satisfies (\ref{the}) and let $\delta_0>0$ be the respective value of $\delta>0$ such that (\ref{the0}) holds. Theorem~\ref{implicit} now applies with $n=3,$ $M=|\beta|mq,$ $\xi=x,$ $\zeta=0,$ $L=\beta,$ and with $q\in(0,1)$ and $\delta_0>0$ selected above.
Thus, by theorem~\ref{implicit} 
there exists $\delta\in(0,\delta_0)$ such that for all $|y|<\delta$ the equation $F(x,y)=0$ has a unique solution $|x(y)|\le m |y|^3.$ Moreover, $x$ is differentiable at $|y|<\delta$.

\vskip0.2cm

\noindent {\it Step 2. Strict monotonicity of $x(y)$.} Differentiating $F(x(y),y)=0$ and expressing $x'(y)$, one concludes
$$
  \lim\limits_{y\to 0}\dfrac{x'(y)}{y^2}=-\dfrac{3\alpha}{\beta}.
$$
Therefore we can assume that $\delta>0$ is so small that $x(y)$ is strictly monotone on $[-\delta,\delta]$. Defining 
$
  \gamma=\min\{|x(-\delta)|,|x(\delta)|\}
$
we conclude that for any $|x|\le\gamma$ the equation $F(x,y)=0$ has a unique solution  $y\in I_{-1}\cup I_{1}$ given by $y=x^{-1}(x).$ Put $y(x)=x^{-1}(x).$

\vskip0.2cm

\noindent {\it Step 3. Local stability of $y(x)$.} Computing the derivative of $P'(y(x))$ we conclude that $y(x)$ is an attractor or repeller according to whether $2\alpha y(x)-\beta \dfrac{x}{y(x)^2}+R'_y(x,y(x))<0$ or $2\alpha y(x)-\beta \dfrac{x}{y(x)^2}+R'_y(x,y(x))>0$.  
Considering $y(x)>0$, we divide these two inequalities by $y(x)$ to obtain
\begin{eqnarray}
{\rm attractor\hskip-0.08cm:} \ \ && 2\alpha -\beta \dfrac{x}{y(x)^3}+\dfrac{R'_y(x,y(x))}{y(x)}<0, \label{attr}\\
 {\rm repeller\hskip-0.08cm:} \ \ && 2\alpha -\beta \dfrac{x}{y(x)^3}+\dfrac{R'_y(x,y(x))}{y(x)}>0.\label{repe}
\end{eqnarray}
Based on (\ref{moreover}), the inequalities (\ref{attr}) and (\ref{repe}) converge to $\alpha<0$ and $\alpha>0$ respectively as $x\to 0$.

\vskip0.1cm

\noindent For those values of $x$ for which $y(x)<0$, the signs in inequalities (\ref{attr}) and (\ref{repe}) flip and $y(x)$ turns out to be an attractor or a repeller according to whether $\alpha>0$ or $\alpha<0.$ 

\vskip0.1cm

\noindent Combining the condition on the sign of $\alpha$ with the condition (\ref{ie}) on the sign of $y(x)$ we conclude that $y(x)$ is an attractor or a repller according to whether $\beta x>0$ or $\beta x<0.$ 

\begin{figure}[h]\label{newfig1}
\begin{center}\includegraphics[scale=0.9]{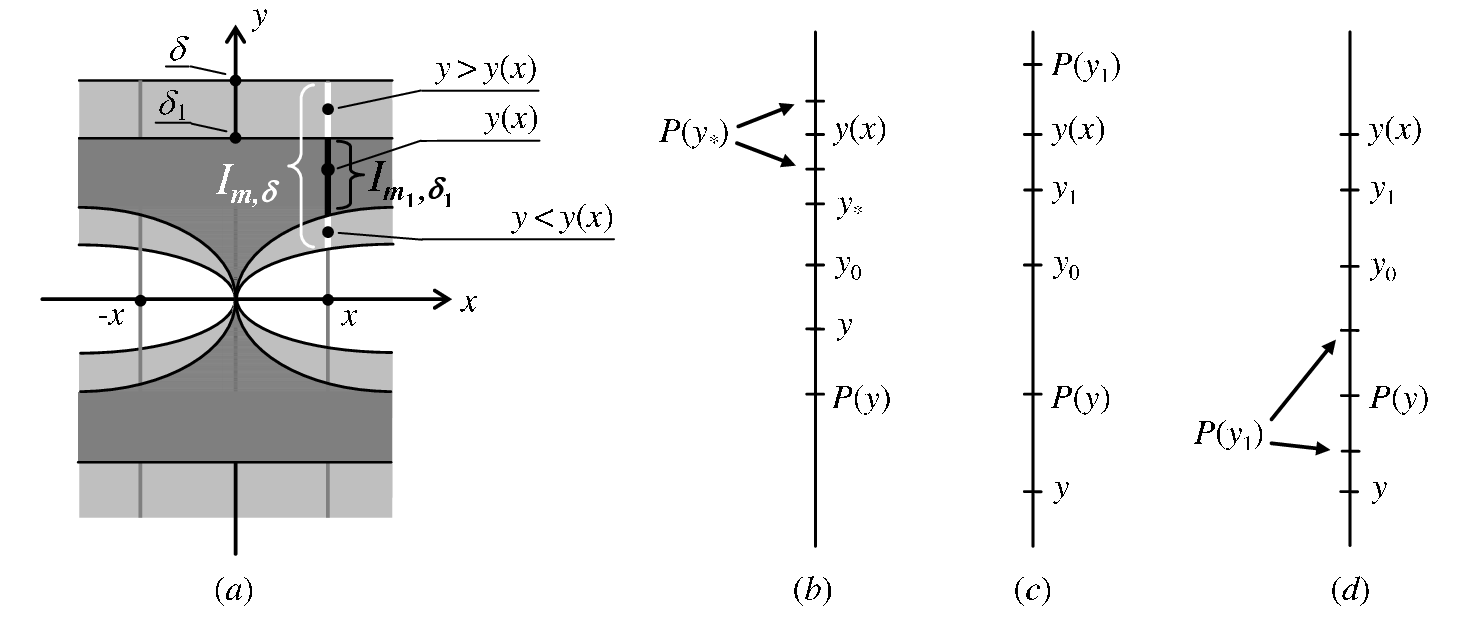}\end{center}
\caption{\footnotesize Illustration of the proof of step 4 of proposition~\ref{dynamicsP}: sets $I_{m,\delta}$ and $I_{m_1,\delta_1}$ (a), step 4.B (b), step 4.C (drawings (c) and (d)).}
\end{figure}

\vskip0.2cm

\noindent {\it Step 4. Attractivity of $y(x)$ in $I_{{\rm sign}(-\alpha\beta x)}$.}   Let $\delta,m>0$ be some constants for which the conclusions of steps 1-3 hold. Let $\delta_1,m_1>0$, $0<\delta_1<\delta$ and $m>m_1>\dfrac{|\alpha|}{|\beta|}$, be another pair of constants for which the conclusions of steps 1-3 hold too. The proof of attractivity is split into 3 parts.
\begin{itemize}





\item[A:] If $y\in I_{m,\delta}(x)\backslash I_{m_1,\delta_1}(x)$ then $y<y(x)\Longrightarrow P(y)<y(x)$ and $y>y(x)\Longrightarrow P(y)>y(x)$, see Fig.~\ref{newfig1}. 

{\it Proof.} Consider $y>0$. Let us show that $P(y)<y(x)$ for any $\sqrt[3]{|x|/m}\le y\le \sqrt[3]{|x|/m_1}$. We have
\begin{eqnarray*}
& &   \hskip-0.5cm\left|\dfrac{P(y)}{y(x)}\right|=\left|\dfrac{y+\alpha y^2+\beta x/y+R(x,y)}{y(x)}\right|\le\\
   && \hskip-0.5cm\le\sqrt[3]{\left|\dfrac{x}{y(x)^3}\right|}\cdot\dfrac{1}{\sqrt[3]{m_1}}\left(1+|\alpha|\sqrt[3]{|x|/m_1}+|\beta|\dfrac{|x|}{\sqrt[3]{|x|^2/m^2}}+\dfrac{|R(x,y)|}{|y|}\cdot\dfrac{\sqrt[3]{|x|/m_1}}{\sqrt[3]{|x|/m}}\right).
\end{eqnarray*}
Since $m_1>|\alpha|/|\beta|,$ then $\sqrt[3]{{x/y(x)^3}}\cdot{\sqrt[3]{1/m_1}}\to k<1$ as $x\to 0.$ Therefore, $P(y)/y(x)<1$, if $|x|\le\gamma$ and $\gamma>0$ is sufficiently small. Let us now $\delta_1\le y\le \delta$ and $k\in(0,1)$ is arbitrary. Then 
$$
P(y)=y+\alpha y^2+ \beta \dfrac{x}{y}+R(x,y)=y\left(1+\alpha y+\beta\dfrac{x}{y^2}+\dfrac{R(x,y)}{y}\right)\ge \delta_1\cdot k
$$
for all $|x|\le m|y|^3,$ $|y|\le\delta,$ provided that $\delta>0$ was chosen sufficiently small. By diminishing $\gamma>0$ we can get $|y(x)|<\delta_1 k,$ so that $P(y)>y(x)$.


\vskip0.1cm

The case when $y<0$ can be considered by analogy. 

\vskip0.1cm

\item[B:] If $y\in I_{m,\delta}(x)\backslash I_{m_1,\delta_1}(x)$ then $|P(y)-y(x)|<|y-y(x)|.$

{\it Proof.} Consider the case where $y(x)>0$. Let $0<y<y(x).$ From part A we have that $P(y)<y(x),$ so it is enough to show that $y<P(y).$ Assume the contrary, i.e. that $P(y)<y$. Step 3 implies that there exists $y<y_*<y(x)$ such that $P(y_*)>y_*$, which implies the existence of a fixed point $y_0\in(y,y_*)$ such that $P(y_0)=y_0$, see Figure~\ref{newfig1}(b). This contradicts the uniqueness of $y(x)$ established in Steps 1-2. The case where $y>y(x)$ and the case where $y(x)<0$ can be considered by analogy. 

\vskip0.1cm

\item[C:] If $y\in I_{m,\delta}(x)$ then $|P(y)-y(x)|<|y-y(x)|.$

{\it Proof.} Let $y\in I_{m,\delta}(x)\backslash I_{m_1,\delta_1}(x)$ and $y<y(x).$ Part B implies that 
\begin{equation}\label{STA}
   y<P(y)<y(x).
\end{equation}
Let us show that (\ref{STA}) implies that 
\begin{equation}\label{COR}
   y_1<P(y_1)<y(x),\quad{\rm for \ all\ }\quad y_1\in (y,y(x)).
\end{equation}
Assume the contrary, i.e. assume that there exists $y_1\in(y,y(x))$ such that $P(y_1)\not\in(y_1,y(x))$.

\vskip0.1cm

Case 1: $P(y_1)>y(x)$, see Figure~\ref{newfig1}(c). In this case there exists $y_0\in(y,y_1)$ such that $P(y_0)=y(x),$ which cannot happen because $y(x)$ is a fixed point of $P.$ 

\vskip0.1cm

Case 2: $P(y_1)<y(x)$, see Figure~\ref{newfig1}(d). This implies the existence of $y_0\in(y,y_1)$ such that $P(y_0)=y_0$ which contradicts the uniqueness of the fixed point $y(x)$ in $I_{m,\delta}.$

\vskip0.1cm

Therefore, (\ref{COR}) is a correct statement. 

\vskip0.1cm

Along the same lines one can prove that if $y(x)<P(y)<y$ then $y(x)<P(y_1)<y_1$ for all $y\in(y(x),y).$

\end{itemize} 

\vskip0.1cm

\noindent Part C implies that $P^n(y)\to y(x)$ as $n\to\infty$ for any $y\in I_{m,\delta}.$

\vskip0.2cm

\noindent The repelling statement can be proved by analogy. \qed

\subsection{The main theorem}\label{secmain}

$ \ $

\vskip0.1cm

We are finally ready to formulate the main result of this paper.

\vskip0.1cm

\begin{theorem}\label{mainthm} Let the $C^4$ functions $f^L,$ $g^L,$ $f^R,$ $f^L$ satisfy (\ref{ff}), i.e. the origin is a fold-fold singularity of the reduced system (\ref{np}).  Assume that 
\begin{eqnarray}
&& \hskip0.4cm f^R{}'_y(0)f^L{}'_y(0)>0,\ \ g^R(0)g^L(0)<0,\label{c2}\\
&& \hskip0.4cm\alpha\not=0,\ \  {\rm where}\ \alpha=\alpha^L-\alpha^R,\ \  \alpha^{L,R}=2\dfrac{f^{L,R}{}'_x(0)+g^{L,R}{}'_y(0)}{g^{L,R}(0)}+\dfrac{f^{L,R}{}'_y{}'_y(0)}{f^{L,R}{}'_y(0)},\label{c3} \\
&& \hskip0.4cm
\alpha\beta f^R{}'_y(0)<0,\ \ {\rm where}\ \beta=\beta^R-\beta^L\ \ {\rm with}\ \ \beta^{L,R}=-\dfrac{2g^{L,R}(0)}{f^{L,R}{}'_y(0)}.\label{c4}
\end{eqnarray}
Then, for any $m>|\alpha|/|\beta|$ there exists $\delta>0$ such that for $xf^R{}'_y(0)g^R(0)<0$ sufficiently close to zero, the switched system 
(\ref{fg})-(\ref{RL}) admits a unique $$y\in J=\left\{y:-\sign(g^R(0))y\in\left[\sqrt[3]{|x/m|},\delta\right]\right\}$$ such that $(x,y)$ is the initial condition of a limit cycle of (\ref{fg})-(\ref{RL}). If $(x(t),y(t))$ is any solution of (\ref{fg})-(\ref{RL}) with 
$(x(0),y(0))\in \{x\}\times J$ then the sequence $\bigcup\limits_{t\ge 0}\left(\{(x(t),y(t))\}\cap \{x\}\times J\right)$ accumulates at $(x,y)$ or contains only a finite number of elements  
according to whether $\alpha g^R(0)>0$ or $\alpha g^R(0)<0.$ In particular, the limit cycle is orbitally stable, if $\alpha g^R(0)>0$,  and unstable, if $\alpha g^R(0)<0.$
\end{theorem}

\vskip0.2cm

{\it Proof.}  The conditions  (\ref{c2})-(\ref{c3}) just restate what is assumed in corollary~\ref{normalP} and proposition~\ref{dynamicsP} already. The condition  (\ref{c4}) implies that $\sign(-\alpha\beta x)$ in (\ref{ie}) with $\sign(x)$ given by (\ref{inadd}) coincides with $-\sign(g^R(0))$ in (\ref{thenfor}), so that the fixed point $y(x)$ of $P$ given by proposition~\ref{dynamicsP} lives in the domain (\ref{thenfor}) where $P$ is the Poincar\'e map of switched system (\ref{fg})-(\ref{RL}) as established in corollary~\ref{normalP}. Furthermore, due to (\ref{c4}), the sign of $x\beta$ with $\sign(x)$ given by (\ref{inadd}) coincides with the sign of $\alpha g^R(0).$ The conclusion, therefore, follows by combining the statements of corollary~\ref{normalP} and proposition~\ref{dynamicsP}.\qed

\vskip0.2cm

\noindent We refer the reader to remark~\ref{note} to clarify that there is no uncertainty in the choice of the initial letter $i$ in (\ref{fg}) when considering the solution $(x(t),y(t))$ in theorem~\ref{mainthm}. This initial letter is always R by construction.

\pagebreak

\section{Applications}\label{appl}

\subsection{A switched mass-spring oscillator}\label{appl1} 
\begin{figure}[h]\label{f3}
\begin{center}\includegraphics[scale=0.9]{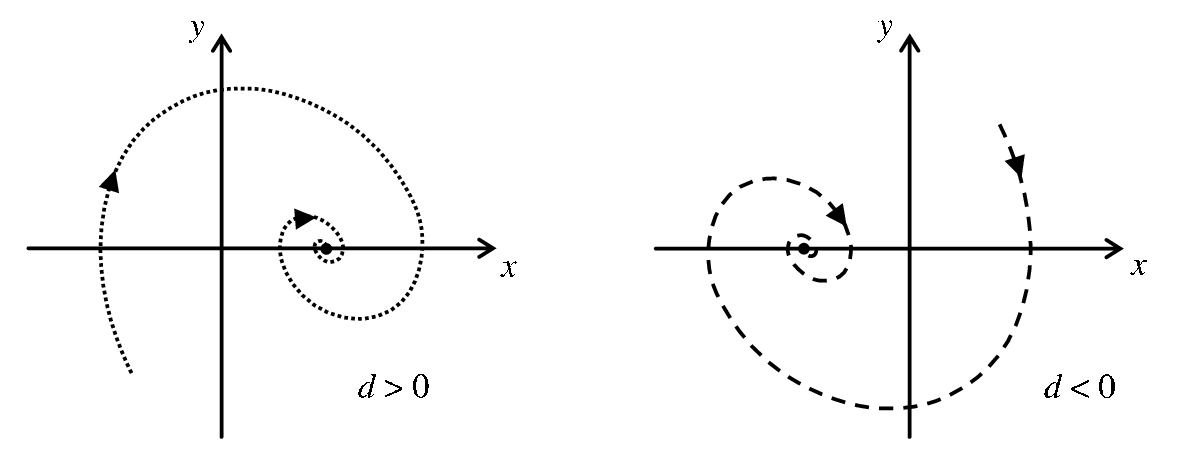}\end{center}
\caption{\footnotesize Sample trajectories of (\ref{mso}) for positive and negative values of $d$.}
\end{figure}

A toy model that illustrates the essence of our main result is the simple mass-spring oscillator 
\begin{equation}\label{mso}
  \begin{array}{l}
      \dot x(t)=y(t),\\
      \dot y(t)=-x(t)-c y(t)+d,
  \end{array}
\end{equation}
whose damping $c$ and forcing $d$
switch from $(c^L,d^L)$ to $(c^R,d^R)$ and back. As a consequence, the unique (globally attracting) equilibrium of this oscillator alternates between $(x,y)=(d^L,0)$ and $(x,y)=(d^R,0),$ see Figure~\ref{f3}. 
We will equip oscillator (\ref{mso}) with the following switching law
\begin{equation}\label{sp}
  \begin{array}{l}
      (c,d):=(c^L,d^L),\quad{\rm if}\ x(t)=-x,\\
      (c,d):=(c^R,d^R),\quad{\rm if}\ x(t)=x.
  \end{array}
\end{equation}
The switching law (\ref{sp}) can be, for example, executed by switching magnets (assuming that the mass in the oscillator (\ref{mso}) is metal) which are connected with sensors positioned at coordinates $-x$ and $x$, see \cite{lamb}.

\begin{proposition}\label{prosp} Consider $c^L>0$, $c^R>0$, $d^L d^R<0$ and introduce
$$
  a=-2\dfrac{c^L}{d^L}+2\dfrac{c^R}{d^R}, \quad b=-2d^R+2d^L.
$$
Then, for any  $m>|a|/|b|$ and for $xd^R<0$ sufficiently close to zero, the switched mass-spring oscillator (\ref{mso})-(\ref{sp}) admits a unique limit cycle with the initial condition $(x,y(x))\to 0$ as $x\to 0$ and $|x|/|y(x)|^3\le m.$ The limit cycle is orbitally stable. 

\end{proposition}

\vskip0.2cm

{\it Proof.} We have
\begin{eqnarray*}\hskip-0.2cm
  \left(\begin{array}{c}
    f^L(0)\\ g^L(0)
\end{array}
\right)= \left(\begin{array}{c}
    0\\ d^L
\end{array}
\right), && \quad  \left(\begin{array}{c}
    f^L{}'(0)\\ g^L{}'(0)
\end{array}
\right)= \left(\begin{array}{cc}
    0 & 1\\ -1 & -c^L
\end{array}
\right),\\ \left(\begin{array}{c}
    f^R(0)\\ g^R(0)
\end{array}
\right)= \left(\begin{array}{c}
    0\\ d^R
\end{array}
\right), && \quad \left(\begin{array}{c}
    f^R{}'(0)\\ g^R{}'(0)
\end{array}
\right)= \left(\begin{array}{cc}
    0 & 1\\ -1 & -c^R
\end{array}
\right).
\end{eqnarray*}
The two inequalities in (\ref{c2}) become $1>0$ and $d^Ld^R<0.$ The constants and the inequality in (\ref{c3}) turn out to be
$$
  \alpha^L=2\dfrac{-c^L}{d^L},\quad   \alpha^R=2\dfrac{-c^R}{d^R},\quad ad^R\not=0.
$$
The constants and the inequality in (\ref{c4}) evaluate to
$$
   \beta^L=-{2d^L},\quad \beta^R=-2d^R,\quad \left(-\dfrac{c^L}{d^L}+\dfrac{c^R}{d^R}\right)\left(-d^R+d^L\right)<0.
$$
The condition for stability leads to $ad^R>0$ or, equivalently, to $(-d^R+d^L)d^R<0,$ which always holds.\qed

\vskip0.2cm

\noindent Proposition~\ref{prosp} supports the simulations carried out for (\ref{sp}) in \cite{lamb}.

\subsection{An anti-lock braking system}\label{appl2}

The analysis in this section concerns  the {\it single-corner}  model (also known as {\it quater-car} model). This model 
is typically used for the preliminary design of wheel braking control algorithms (see \cite{pas06,tan09}). 
If the longitudinal dynamics of the vehicle  is much slower than the rotational dynamics of the wheel than the interplay between the longitudinal slip and the braking torque of the wheel reads as 
(see \cite{tan09}) 
\begin{equation}\label{abs}
   \begin{array}{rcl}
      \dot\lambda(t)&=&-F(\lambda(t))+F_0 T_b(t),
            \end{array}
\end{equation}
with
$$
   F(\lambda)=\dfrac{1}{\nu}\left(\dfrac{1-\lambda}{m}+\dfrac{r^2}{J}\right)F_z\mu(\lambda),\quad \mu(\lambda)=\theta_{r1}(1-{\rm exp}(-\lambda\theta_{r2}))-\lambda\theta_{r3},\quad F_0=\dfrac{r}{\nu J},
$$
where $\nu$  is the longitudinal speed of the vehicle,
$r$ is the wheel radius,
$J$ is the moment of inertia of the wheel,
$m$ is the mass of the quarter-car,
$F_z$ is the vertical force at the tire-road contact point, $\theta_{r1},$ $\theta_{r2},$ $\theta_{r3}$ are positive constants that reflect the road conditions.

\vskip0.2cm

\begin{figure}[h]\label{f4}
\begin{center}\includegraphics[scale=0.9]{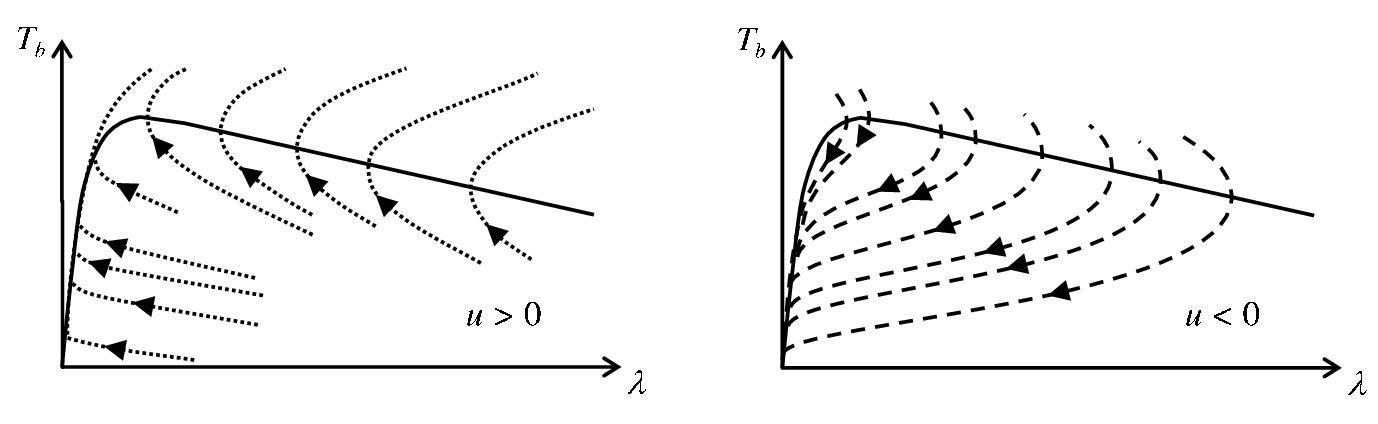}\end{center}
\caption{\footnotesize Trajectories of system (\ref{abs})-(\ref{absa}) for $k<0$ and $k>0$. The solid line is the graph of the function $T_b=\dfrac{1}{F_0}F(\lambda)$, i.e. the set of $(\lambda,T_b)$ where the first component of the right-hand-side of (\ref{abs})-(\ref{absa}) vanishes.}
\end{figure}

\noindent The assumption of either a fixed or a user-selectable actuator rate limit has been used in several works on rule-based ABS, see 
\cite{pas06,tan09} and references therein. This leads to the following differential equation for $T_b$
\begin{equation}\tag{\ref{abs}a}\label{absa}
  \dot T_b=u,
\end{equation}
 where $u$ is a control input, that can take a finite number of values. In this paper we consider the simplest control logic where the actuator rates can take two values $k_1 = k$ and $k_2 = -k$ of equal modulus, see Figure~\ref{f4} for the phase portrait of system (\ref{abs})-(\ref{absa}).  Following \cite{pas06,tan09}, the goal of the controller is to make $\lambda\mapsto\lambda(t)$ oscillating periodically within the interval $[\lambda_0-\Delta\lambda,\lambda_0+\Delta\lambda]$, where $\lambda_0$ is the measured wheel slip and $\Delta\lambda$ is the measurement error, i.e. the deviation of the measured wheel slip from the optimal value (that ensures the fastest braking).
Whenever $\lambda(t)$ hits $\lambda_0+\Delta\lambda$, the change of the braking torque $\dot{T}_b = u$ shall switch from $u= k_1>0$ to $u= k_2<0$. On the contrary, as soon as $\lambda_0-\Delta\lambda$ is hit, it shall switch from $u = k_2$ back to $u=k_1$. This control logic will be shortly formulated as
\begin{equation}\label{absc}
   \begin{array}{rcl}
      \dot T_b(t)&:=&\left\{\begin{array}{ll}
              -k, &{\rm if}\ \lambda(t)\ge \lambda_0+\Delta\lambda,\\
             k, &{\rm if}\ \lambda(t)\le \lambda_0-\Delta\lambda.
            \end{array}\right.
   \end{array}
\end{equation}

\noindent The existence of limit cycles in the anti-lock braking system (\ref{abs})-(\ref{absc}) has been recently established in \cite{koppen} in the case where the magnitude of the right-hand-side in (\ref{abs}) is small, but the $\Delta \lambda$ is fixed. In contrast, our next proposition deals with arbitrary right-hand-sides of (\ref{abs}), but needs smallness of  $\Delta \lambda$.

\begin{proposition}\label{propabs} Assume that $F_0k\not=0$. Fix $\lambda_0\in\mathbb{R}$ and consider
$$
   a=4\dfrac{F'(\lambda_0)}{k},\quad b=-4\dfrac{k}{F_0}.
$$
If $F'(\lambda_0)>0$ then for any  $m>|a|/|b|$ and for $\Delta\lambda\cdot k>0$ sufficiently close to zero, the anti-lock braking system (\ref{abs})-(\ref{absc}) admits a unique limit cycle with the initial condition $(x,y(x))\to 0$ as $x\to 0$ and $|x|/|y(x)|^3\le m.$ The limit cycle is orbitally stable. 
\end{proposition}

\vskip0.2cm
{\it Proof.} The first component of the right-hand-side vanishes when $T_b=(1/F_0)F(\lambda).$ Therefore,  points of the type $(\lambda,F_0^{-1}F(\lambda))$ are potential candidates to produce limit cycles of (\ref{abs}). We have 
$$
\begin{array}{ll}
  \left(\begin{array}{c}
    f^L(\lambda,F_0^{-1}F(\lambda))\\ g^L(\lambda,F_0^{-1}F(\lambda))
\end{array}
\right)= \left(\begin{array}{c}
    0\\ -k
\end{array}
\right),\  &  \left(\begin{array}{c}
    f^L{}'(\lambda,F_0^{-1}F(\lambda))\\ g^L{}'(\lambda,F_0^{-1}F(\lambda))
\end{array}
\right)= \left(\begin{array}{cc}
    -F'(\lambda) & F_0\\ 0 & 0
\end{array}
\right),\\ 
\left(\begin{array}{c}
    f^R(\lambda,F_0^{-1}F(\lambda))\\ g^R(\lambda,F_0^{-1}F(\lambda))
\end{array}
\right)= \left(\begin{array}{c}
    0\\ k
\end{array}
\right),\ &  \left(\begin{array}{c}
    f^R{}'(\lambda,F_0^{-1}F(\lambda))\\ g^R{}'(\lambda,F_0^{-1}F(\lambda))
\end{array}
\right)= \left(\begin{array}{cc}
    -F'(\lambda) & F_0\\ 0 & 0
\end{array}
\right).
\end{array}
$$
The two inequalities in (\ref{c2}) become $(F_0)^2>0$ and $-k^2<0.$ The constants and the inequality in (\ref{c3}) turn out to be
$$
  \alpha^L=2\dfrac{-F'(\lambda)}{k},\quad   \alpha^R=2\dfrac{-F'(\lambda)}{-k},\quad a=-4\dfrac{F'(\lambda)}{k},\quad F'(\lambda)\not=0.
$$
The constants and the inequality in (\ref{c4}) evaluate to
$$
   \beta^L=-\dfrac{2k}{F_0},\quad \beta^R=-\dfrac{-2k}{F_0},\quad b=\dfrac{4k}{F_0},\quad -F'(\lambda)<0,
$$
which implies that the stability condition holds true.
\qed

\section{Conclusions} \noindent The main result of this paper can be further extended to a greater number of vector fields and switching thresholds that govern the switched system (\ref{fg})-(\ref{RL}). Building upon the ideas of section~\ref{sec2}, such an extension will be just technical as long as planar vector fields are concerned. In particular, the ABS application of section~\ref{appl2} can be extended to controllers with a greater number of discrete actions. A challenging next step in our research is to discover an extension of  the main result to 3-dimensional systems  (\ref{fg})-(\ref{RL}), where the switching lines $\{x\}\times \R$ and $\{-x\}\times \R$ are parallel switching hyperplanes that approach each other and merge into a single hyperplane (say, $L$) as  $x\to 0$, see \cite[Figure~10]{lamb}. The respective reduced system (\ref{np}) will now have a so-called $U$-singularity at the origin \cite{fil} (also known as Teixeira-singularity \cite{teixeira}), where the vector fields on the two sides of $L$ are parallel to $L$. Understanding the dynamics of (\ref{fg})-(\ref{RL}) in 3D will thus need to be built upon the knowledge of the dynamics of $U$-singularity gained since \cite{teixeira1,fil}. The prospective to explore the interesting geometry of $U$-singularity was another reason for us to focus the analysis on the fold-fold singularity, i.e. to work in the setting (\ref{ff}). Moreover, the only part of our proofs which is intrinsically 2-dimensional is the extension of local stability to attractivity in a larger set (step~4 of the proof of proposition~\ref{dynamicsP}).

\appendix\section {Partial derivatives of the general solution at a fold singularity}\label{partial}

Let $t\mapsto (X(t,x,y),Y(t,x,y))$ be the general solution of (\ref{neim}) that satisfies $f(0)=0.$ Then, by direct computation, \hskip2cm\begin{table}[h] 
\def\arraystretch{1.5}
\hskip0.4cm\begin{tabular}{lll}
$X{}_t\hskip-0.1cm'\hskip0.05cm(0)=0,$   & $X{}'_t{}'_x(0)=f{}_x'(0),$ & \hskip-2cm$X{}'_t{}'_y{}'_y(0)=f{}'_y{}'_y(0),$ 
\\
$X{}_x\hskip-0.1cm'\hskip0.05cm(0)=1,$  &  $X{}'_t{}'_y(0,x,y)=f{}'_y(x,y),$ & \hskip-2cm$X{}'_t{}'_t{}'_y(0)=f{}'_x(0)f{}'_y(0)+
$
\\
$X{}_y\hskip-0.1cm'\hskip0.05cm(0)=0,$   & $Y{}'_t{}'_x(0)=g{}'_x(0),$ & \hskip-2cm\hskip0.2cm$+f{}'_y{}'_y(0)g(0)+f{}'_y(0) g{}'_y(0),$ 
\\
$Y{}_t\hskip-0.1cm'\hskip0.05cm(0)=g(0),$  &  $Y{}'_t{}'_y(0)=g{}'_y(0),$ & \hskip-2cm $X{}'_t{}'_t{}'_t(0)=\left[f{}'_x(0)f{}'_y(0)+\right.$
\\
$Y{}_x\hskip-0.1cm'\hskip0.05cm(0)=0,$   & $X{}'_t{}'_t(0,x,y)=f{}'_y(x,y)g(x,y),$ & \hskip-2cm\hskip0.2cm$\left.+
f{}'_y{}'_y(0)g(0)+f{}'_y(0) g{}'_y(0)\right]g(0).$
\\
$Y{}_y\hskip-0.1cm'\hskip0.05cm(0)=1,$ & $Y{}'_t{}'_t(0)=g{}'_y(0)g(0),$ & 
\\
 & $X{}'_x{}'_x(0)=X{}'_y{}'_y(0)=Y{}'_x{}'_x(0)=Y{}'_y{}'_y(0)=0$ & 
\end{tabular}
\end{table}

\section{Limiting relations for the reminders}

\subsection{The reminder $\rho$ in the proof of lemma~\ref{mathcalP}}\label{rhoprop} $ $

{\bf 1.} By construction, $\rho(t,x,y)=\dfrac{X(t,x,y)-x}{t}-X'_t(0,x,y)-X{}'_t{}'_t\dfrac{t}{2}.$ 

\vskip0.2cm

Therefore, by l'Hospital rule,

$\lim\limits_{t\to 0}\rho(t,x,y)=X'_t(0,x,y)-X'(0,x,y)-X{}'_t{}'_t(0,x,y)\cdot 0=0.$

\vskip0.2cm

{\bf 2.} Since $\rho'_t(t,x,y)=\left(\dfrac{X(t,x,y)-x}{t}\right)'-X{}'_t{}'_t(0,x,y)\dfrac{1}{2}=$

$=\dfrac{X'_t(t,x,y)t-(X(t,x,y)-x)}{t^2}-X{}'_t{}'_t(0,x,y)\dfrac{1}{2}$, the l'Hospital rule yields

$\lim\limits_{t\to 0}\rho'_t(t,x,y)=\lim\limits_{t\to 0}\dfrac{X{}'_t{}'_t(t,x,y)t+X'_t(t,x,y)-X'_t(t,x,y)}{2t}-X{}'_t{}'_t(0,x,y)\dfrac{1}{2}=0.$

\vskip0.2cm

{\bf 3.} $\lim\limits_{t\to 0}\rho'_y(t,x,y)=\lim\limits_{t\to 0}\left(\dfrac{X'_y(t,x,y)}{t}-X{}'_t{}'_y(0,x,y)-X{}'_t{}'_y(0,x,y)\dfrac{t}{2}\right)=0.$

\vskip0.2cm

{\bf 4.} One has $\rho{}'_t{}'_y(t,x,y)=\dfrac{X{}'_t{}'_y(t,x,y)t-X'_y(t,x,y)}{t^2}-X{}'_t{}'_t{}'_y(0,x,y)\dfrac{1}{2}$ and by the l'Hospital rule

\noindent $\lim\limits_{t\to 0}\rho{}'_t{}'_y(t,x,y)=\lim\limits_{t\to 0}\dfrac{X{}'_t{}'_t{}'_y(t,x,y)t+X{}'_t{}'_y(t,x,y)-X{}'_t{}'_y(t,x,y)}{2t}-X{}'_t{}'_t{}'_y(0,x,y)\dfrac{1}{2}=0.$

\vskip0.2cm

{\bf 5.}  $\lim\limits_{t\to 0}\rho{}'_y{}'_y(t,x,y)=\lim\limits_{t\to 0}\left(\dfrac{X{}'_t{}'_y(t,x,y)}{t}-X{}'_t{}'_y{}'_y(0,x,y)-X{}'_t{}'_t{}'_y{}'_y(0,x,y)\dfrac{t}{2}\right)=0.$

\vskip0.2cm

{\bf 6.} $\lim\limits_{t\to 0}\rho{}'_x(t,x,y)=\lim\limits_{t\to 0}\left(\dfrac{X'_x(t,x,y)-1}{t}-X{}'_t{}'_x(0,x,y)-X{}'_t{}'_t{}'_x(0,x,y)\dfrac{t}{2}\right)=0.$

\vskip0.2cm

{\bf 7.} $\lim\limits_{t\to 0}\rho{}'_t{}'_x(t,x,y)=\lim\limits_{t\to 0}\left(\dfrac{X{}'_t{}'_x(t,x,y)t-X'_x(t,x,y)+1}{t^2}-X{}'_t{}'_t{}'_x(0,x,y)\dfrac{1}{2}\right)=$

$=\lim\limits_{t\to 0}\dfrac{X{}'_t{}'_t{}'_x(t,x,y)t+X{}'_t{}'_x(t,x,y)-X{}'_t{}'_x(t,x,y)}{2t}-X{}'_t{}'_t{}'_x(0,x,y)\dfrac{1}{2}=0.$

\vskip0.2cm

{\bf 8.} $\lim\limits_{t\to 0}\rho{}'_x{}'_x(t,x,y)=\lim\limits_{t\to 0}\left(\dfrac{X{}'_x{}'_x(t,x,y)}{t}-X{}'_t{}'_x{}'_x(0,x,y)-X{}'_t{}'_t{}'_x{}'_x(0,x,y)\dfrac{t}{2}\right)=0.$

\vskip0.2cm

{\bf 9.} $\lim\limits_{t\to 0}\rho{}'_t{}'_t(t,x,y)=$

\noindent $=\lim\limits_{t\to 0}\dfrac{(X{}'_t{}'_t(t,x,y)t+X'_t(t,x,y)-X'_t(t,x,y))t^2-2t(X'_t(t,x,y)t-(X(t,x,y)-x))}{t^4}=$

\noindent $=\lim\limits_{t\to 0}\dfrac{X{}'_t{}'_t{}'_t(t,x,y)t^3+X{}'_t{}'_t(t,x,y)t^2-2tX{}'_t(t,x,y)+2(X(t,x,y)-x)}{4t^3}=$

\noindent$=\lim\limits_{t\to 0}\dfrac{X{}'_t{}'_t{}'_t{}'_t(t,x,y)t^3+X{}'_t{}'_t{}'_t(t,x,y)4t^2}{12t^2}.$

\subsection{The reminder $\widetilde{\mathcal{R}}$ in the proof of lemma~\ref{lem2}}\label{appendixC} $ $

In contrast to the reminder $\rho$ from lemma~\ref{mathcalP},  the reminder $\widetilde{\mathcal{R}}$ doesn't come from a formal Taylor expansion. That is why we didn't manage to prove limiting properties of $\rho$ and limiting properties of $\widetilde{\mathcal{R}}$ simultaneously. That is why an independent proof for $\widetilde{\mathcal{R}}$ appears below.

{\bf 1.} $\dfrac{\widetilde {\mathcal{R}}(x,y)}{y^2}=\left(-\dfrac{2Y'_t(t_*,x,y)}{\dfrac{X'_t(0,x,y)}{y}}+\dfrac{2Y'_t(0)}{X{}'_t{}'_y(0)}\right)\dfrac{x}{y^3}-\dfrac{Y'_t(t_*,x,y)}{\dfrac{X'_t(0,x,y)}{y}}\cdot\dfrac{\widetilde \Delta(\widetilde T(x,y),x,y)}{y^3}.
$

The later expression approaches zero as $|x|\le m |y|^3$, $y\to 0$, because
\begin{equation}\label{Deltay3}
\begin{array}{l}
\left|\dfrac{\widetilde \Delta(\widetilde T(x,y),x,y)}{y^3}\right|=\left|\dfrac{X(\widetilde T(x,y),x,y)-x-X'_t(0,x,y)\widetilde T(x,y)}{y^3}\right|=\\
=\left|\dfrac{X{}'_t{}'_t(t_{**},x,y)t_* \widetilde T(x,y)}{y^3}\right|\le   |X{}'_t{}'_t(t_{**},x,y)|\cdot\dfrac{(4m)^2}{|f'_y(0)|^2}\cdot |y|.
\end{array}
\end{equation}

\vskip0.2cm

 {\bf 2.}  Differentiating (\ref{2x}) with respect to $y$, 
\begin{eqnarray*}
    \dfrac{\widetilde{T}'_y(x,y)}{y}&=&\dfrac{1}{D}\left(-X{}'_t{}'_y(0,x,y)\dfrac{\widetilde{T}(x,y)}{y^2}-\dfrac{\widetilde\Delta'_y(\widetilde T(x,y),x,y)}{y^2}\right),\\
&& {\rm where }\ \ D=\dfrac{X{}'_t(0,x,y)}{y}+\dfrac{\widetilde\Delta'_t(\widetilde T(x,y),x,y)}{y}
\end{eqnarray*}
Differentiating (\ref{thfor}) with respect to $y$,
\begin{eqnarray*}
\dfrac{\widetilde{\mathcal{R}}'_y(x,y)}{y}&=&\dfrac{Y'_y(\widetilde T(x,y),x,y)-1}{y}+Y'_t(\widetilde T(x,y),x,y)\dfrac{\widetilde{T}'_y(x,y)}{y}-\dfrac{2Y'_t(0)}{X{}'_t{}'_y(0)}\cdot \dfrac{x}{y^3}=\\
&=&Y{}'_y{}'_y(t_*,x,y)\dfrac{\widetilde T(x,y)}{y}+\left(\dfrac{Y'_t(\widetilde T(x,y),x,y)}{D}\cdot\dfrac{2X{}'_t{}'_y(0,x,y)}{\dfrac{X'_t(0,x,y)}{y}}-\dfrac{2Y'_t(0)}{X{}'_t{}'_y(0)}\right)\dfrac{x}{y^3} +\\
&&+\dfrac{1}{D}\left(X{}'_t{}'_y(0,x,y)\dfrac{\widetilde \Delta(\widetilde{T}(x,y),x,y)}{y^3}\cdot\dfrac{1}{\dfrac{X'_t(0,x,y)}{y}}-\dfrac{\widetilde\Delta'_y(\widetilde T(x,y),x,y)}{y^2}\right).
\end{eqnarray*}

\noindent In what follows we show that 

\centerline{$\lim\limits_{|x|\le m|y|^3,\ y\to 0}{\widetilde\Delta'_t(\widetilde T(x,y),x,y)}/{y}=0$ \quad and\quad $\lim\limits_{|x|\le m|y|^3,\ y\to 0}{\widetilde\Delta'_y(\widetilde T(x,y),x,y)}/{y^2}=0$,} 

\noindent which combined with (\ref{uniquesolution}) and (\ref{Deltay3})
implies that $\lim\limits_{|x|\le m|y|^3,\ y\to 0}\widetilde{\mathcal{R}}'_y(x,y)/y=0.$

\vskip0.2cm

{\bf 3.} Taking into account (\ref{uniquesolution}) one gets

\noindent  $\lim\limits_{|x|\le m|y|^3,\ y\to 0}\dfrac{\widetilde\Delta'_t(\widetilde T(x,y),x,y)}{y}=\lim\limits_{|x|\le m|y|^3,\ y\to 0}\dfrac{X'_t(\widetilde T(x,y),x,y)-X'_t(0,x,y)}{y}=0$.

\vskip0.2cm

{\bf 4.}  $\lim\limits_{|x|\le m|y|^3,\ y\to 0}\dfrac{\widetilde\Delta'_y(\widetilde T(x,y),x,y)}{y^2}=0$ follows by observing that $\dfrac{\widetilde\Delta'_y(t,x,y)}{y^2}=$

$=\dfrac{X'_y(t,x,y)-X{}'_t{}'_y(0,x,y)t}{y^2}=\dfrac{X{}'_t{}'_y(t_*,x,y)t-X{}'_t{}'_y(0,x,y)t}{y^2}=\dfrac{X{}'_t{}'_t{}'_y(t_{**},x,y)t_* t}{y^2}$

and by using (\ref{uniquesolution}) again.

\subsection{The reminder $\mathcal{R}$ in the proof of  corollary~\ref{cor2_3}}\label{reminderR} The limiting relations (\ref{SS}) follow by writing 
\begin{eqnarray*}
 \dfrac{T(x,y)}{y}&=&\dfrac{T'_x(x_*,y)x+T(0,y)}{y},\\
   \dfrac{{\mathcal R}'_y(x,y)}{y}&=&\dfrac{{\mathcal R}'_y(x,y)-{\mathcal R}'_y(x,0)+{\mathcal R}'_y(x,0)}{y}=\mathcal{R}{}'_y{}'_y(x,y_*)+\dfrac{\mathcal{R}{}'_x{}'_y(x_*,0)x}{y},\\
\dfrac{{\mathcal R}(x,y)}{y^2}&=&\dfrac{{\mathcal R}(x,y)-{\mathcal R}(x,0)+{\mathcal R}(x,0)}{y^2}=\dfrac{{\mathcal R}'_y(x,y_*)}{y}+\dfrac{{\mathcal R}(x,0)}{y^2}=\\
&&=\dfrac{{\mathcal R}'_y(x,y_*)}{y}+\dfrac{{\mathcal R}'_x(x_*,0)x}{y^2}.
\end{eqnarray*}

\subsection{The reminder $r$ in the proof of corollary~\ref{compo}}\label{reminderr} The proof of (\ref{limrel}) is split in several steps.

{\bf 1.}  We have
\begin{equation}\label{IN}
\begin{array}{l}
\lim\limits_{|x|\le m|y|^3,\ y\to 0}\widetilde{\mathcal{R}}(x,-y+\alpha y^2+\mathcal{R}(x,y))=\\
=\lim\limits_{|x|\le m\left|\dfrac{y}{-y+\alpha y^2+\mathcal{R}(x,y)}\right|\cdot\left|-y+\alpha y^2+\mathcal{R}(x,y)\right|,\ y\to 0}\hskip-0.5cm\widetilde{\mathcal{R}}(x,-y+\alpha y^2+\mathcal{R}(x,y)).
\end{array}
\end{equation}
Take $M>m.$ Then 
\begin{equation}\label{IN1}
m\left|\dfrac{y}{-y-\alpha y^2+\mathcal{R}(x,y)}\right|\le M\quad{\rm for}\ |x|\le m|y|^3\ {\rm and}\ |y|>0\ {\rm sufficiently\ small.}
\end{equation}
Since
$$
  \lim\limits_{|x|\le M|-y+\alpha y^2+\mathcal{R}(x,y)|,\ -y+\alpha y^2+\mathcal{R}(x,y)\to 0}\widetilde{\mathcal{R}}(x,-y+\alpha y^2+\mathcal{R}(x,y))=0,
$$
then the limit in (\ref{IN}) is zero by (\ref{IN1}).
\

\begin{eqnarray*}
{\bf 2.}\qquad&& \lim\limits_{|x|\le m|y|^3,\ y\to 0}\dfrac{1}{y}\dfrac{d}{dy}\left[\widetilde{\mathcal{R}}(x,-y+\alpha y^2+\mathcal{R}(x,y))\right]=\\
   && =\lim\limits_{|x|\le m|y|^3,\ y\to 0}\dfrac{\widetilde{\mathcal{R}}'_y(x,-y+\alpha y^2+\mathcal{R}(x,y))}{y}(-1+2\alpha y+\mathcal{R}'_y(x,y))=\\
&&  =\lim\limits_{|x|\le m|y|^3,\ y\to 0} \dfrac{\widetilde{\mathcal{R}}'_y(x,-y+\alpha y^2+\mathcal{R}(x,y))}{-y+\alpha y^2+\mathcal{R}(x,y)}\cdot\dfrac{-y+\alpha y^2+\mathcal{R}(x,y)}{y}\cdot \\
&& \qquad \qquad\qquad\cdot(-1+2\alpha y+\mathcal{R}'_y(x,y))=0\cdot (-1)\cdot(-1)
\end{eqnarray*}
by (\ref{tildeR}) and (\ref{SS}).

\begin{eqnarray*}
&&\qquad {\bf 3.}\ \lim\limits_{|x|\le m|y|^3,\ y\to 0}\dfrac{d}{dx}\left[\widetilde{\mathcal{R}}(x,-y+\alpha y^2+\mathcal{R}(x,y))\right]=\\
&& =\lim\limits_{|x|\le m|y|^3,\ y\to 0}\left(\widetilde{\mathcal{R}}'_x(x,-y+\alpha y^2+\mathcal{R}(x,y))+\mathcal{R}'_y(x,-y+\alpha y^2+\mathcal{R}(x,y))\mathcal{R}'_x(x,y)\right),
\end{eqnarray*}
which equals zero by the same reasons as in 2.

$$
\hskip-2.7cm\begin{array}{rcl}
  {\bf 4.} && \lim\limits_{|x|\le m|y|^3,\ y\to 0}\left(\beta\dfrac{x}{y}+\beta\dfrac{x}{-y+\alpha y^2 +\mathcal{R}(x,y)}\right)=\\
&& =\lim\limits_{|x|\le m|y|^3,\ y\to 0}\beta\dfrac{x}{y}\left(1+\dfrac{1}{-1+\alpha y+\dfrac{\mathcal{R}(x,y)}{y}}\right)=0.
\end{array}
$$

$$
\hskip0.3cm\begin{array}{rcl}
{\bf 5.} & &  \lim\limits_{|x|\le m|y|^3,\ y\to 0}\dfrac{1}{y}\dfrac{d}{dy}\left[\beta \dfrac{x}{y}+\beta\dfrac{x}{-y+\alpha y^2+\mathcal{R}(x,y)}\right]=\\
&&\hskip-1.5cm =\lim\limits_{|x|\le m|y|^3,\ y\to 0}\dfrac{x}{y^3}\left(-\beta-\beta\dfrac{1}{\left(-1+\alpha y+\dfrac{\mathcal{R}(x,y)}{y}\right)^2}(-1+2\alpha y+\mathcal{R}'_y(x,y))\right)=0
\end{array}
$$
because $\dfrac{x}{y^3}$ is bounded and the term in the brackets approaches zero as $|x|\le m|y|^3$, $y\to 0.$

$$
\hskip-0.8cm\begin{array}{rcl}
 & &  \hskip0.6cm {\bf 6.} \ \ \lim\limits_{|x|\le m|y|^3,\ y\to 0}\dfrac{d}{dx}\left[\beta\dfrac{x}{y}+\beta\dfrac{x}{-y+\alpha y^2+\mathcal{R}(x,y)}\right]y=\\
&& =\lim\limits_{|x|\le m|y|^3,\ y\to 0}\left(\dfrac{\beta}{y}+\dfrac{\beta}{-y+\alpha y^2+\mathcal{R}(x,y)}-\dfrac{x}{(-y+\alpha y^2+\mathcal{R}(x,y))^2}\mathcal{R}'_x(x,y)\right)y=0.
\end{array}
$$

\section{An implicit function theorem for implicit functions that branch from the boundary of a set}

\vskip0.2cm

\noindent The following technical modification of the standard implicit function theorem is required for the proofs in this paper. This theorem gives a sufficient condition for the existence of an implicit function that branches from the origin in the case where the derivative of the governing equation is not defined at the origin, but is given in a region whose boundary passes through the origin, as in the right drawing of Figure~\ref{fig2}.

\vskip0.2cm

\begin{theorem}\label{implicit}  Consider $F\in C^1(\mathbb{R}\times\mathbb{R}\times \mathbb{R},\mathbb{R})$ such that $F(0,0,0)=0.$ 
Assume that there exist  $M>0$, $m>0,$ and $\delta_0>0$ such that 
\begin{equation}\label{M} 
|F(0,\zeta,y)|\le My^n,\qquad  {\rm for\ any\ } \ |\zeta|\le m |y^3|, \ |y|\le \delta_0,
\end{equation}
where $n\in\mathbb{N}$ is a constant.
If there exists $q\in(0,1)$ such that
\begin{equation}\label{Lq} 
F'_\xi(\xi,\zeta,y)\to L\not=0,\quad {\rm as\ } \ |\xi|\le\dfrac{M}{|L|q}|y|^n,\ |\zeta|\le m|y^3|, \ y\to 0,
\end{equation}
then, for some  $\delta\in(0,\delta_0)$, the equation
\begin{equation}\label{eqF}
   F(\xi,\zeta,y)=0
\end{equation}
has a unique solution $|\xi(\zeta,y)|\le \dfrac{M}{Lq}|y|^n$ in the set 
\begin{equation}\label{set}
|\zeta|\le m|y^3|, \ \ |y|\le \delta.
\end{equation} Moreover, the function $\xi$ is differentiable in the interior of this region.
\end{theorem}

\vskip0.2cm

{\it Proof.} {\it Step 1. The existence.} Consider
$$
  A_{(\zeta,y)}(\xi)=\xi-\dfrac{1}{L}F(\xi,\zeta,y).
$$
We want to show that for any $\delta>0$ sufficiently small the map $A_{(\zeta,y)}$ maps the disk $|\xi|\le\dfrac{M}{|L|q}\delta^n$ into itself and contracts in this disk. We have
$$
  A_{(\zeta,y)}{}'(\xi)=\dfrac{1}{L}\left(L-F'_\xi(\xi,\zeta,y)\right).
$$
Diminishing $\delta>0$ so that
$$
  \left|A_{(\zeta,y)}{}'(\xi)\right|\le 1-q,\quad {\rm for\ any\ }|\xi|\le\dfrac{M}{Lq}|y|^n,\ |\xi|\le m|y^3|,\ |y|\le\delta,
$$
one gets 
$$
  \left|A_{(\zeta,y)}(\xi)\right|\le \left|A_{(\zeta,y)}(0)\right|+\left|A_{(\zeta,y)}(\xi)-A_{(\zeta,y)}(0)\right|\le \dfrac{1}{L}M\delta^n+(1-q)\dfrac{M}{Lq}\delta^n=\dfrac{M}{Lq}\delta^n.
$$

\noindent {\it Step 2. The differentiability.} Let $(\zeta_0,y_0)$ be a point of the interior of (\ref{set}). The equation (\ref{eqF}) can be solved in $\xi$ near $(\zeta_0,y_0)$ because the conditions of the standard implicit function theorem (see e.g. \cite[\S~8.5.4, Theorem~1]{zorich}) hold at $(\xi(\zeta_0,y_0),\zeta_0,y_0)$. The differentiability of $\xi$ at $(\zeta_0,y_0)$, therefore, follows from the classical result on the differentiability of the implicit function (see same source \cite[\S~8.5.4, Theorem~1]{zorich}).

\qed

{\bf Acknowledgments.} This research started during the stay of the author in the University of Cologne under a
Humboldt Postdoctoral Fellowship.
 Different stages of this work were presented at the International Conference on Dynamical Systems celebrating the 70th birthday of Marco Antonio Teixeira (May 26-29, 2014, \href{http://www.ime.unicamp.br/~rmiranda/mat70/MAT70/Welcome.html}{\underline {link}}), in Applied Mathematics Seminar at the University of Colorado at Boulder (June 16, 2014, \href{http://amath.colorado.edu/content/special-oleg-makarenkov}{\underline{link}}), at the International Conference on Differential Equations and Dynamical Systems in Suzdal (July 4-9, 2014, \href{http://agora.guru.ru/display.php?conf=diff2014}{\underline{link}}), in Dynamical Systems Seminar at Imperial College London (Sep 16, 2014, \href{http://wwwf.imperial.ac.uk/~mrasmuss/DynamIC/seminars.php}{\underline{link}}), at  SIAM Conference on Control and its Applications (June 8-10, 2015, \href{http://meetings.siam.org/sess/dsp_talk.cfm?p=71255}{\underline{link}}),
in Dynamical Systems Seminar at the University of Houston (Feb 16, 2016, \href{http://www.math.uh.edu/dynamics/talk.php?file=speaker2016.02.15}{\underline{link}}). The author thanks the organizers for hospitality.

\end{document}